\documentclass[a4paper,12pt
]{article}
\usepackage{amsmath}
\usepackage{amssymb}
\usepackage[dvipdfmx]{graphicx}
\usepackage{graphicx}
\usepackage{here}
\usepackage{mathpazo}
\newtheorem{thm}{Theorem}[section]
\newtheorem{lem}[thm]{Lemma}
\newtheorem{prop}[thm]{Propoition}

\newtheorem{cor}[thm]{Corollary}
\newtheorem{rem}[thm]{Remark}
\newtheorem{df}[thm]{Definition}

\newtheorem{prob}[thm]{Problem}

\newtheorem{constr}[thm]{Construction}

\newcommand{\R}{{\Bbb R}}
\newcommand{\N}{{\Bbb N}}
\newcommand{\C}{{\Bbb C}}
\newcommand{\Z}{{\Bbb Z}}
\newcommand{\Q}{{\Bbb Q}}

\newcommand{\Diff}{{\mathit{Diff}}}
\newcommand{\F}{{\cal F}}
\newcommand{\G}{{\cal G}}

\newcommand{\K}{{\cal K}}

\newcommand{\IE}{{\it i.e.}, }
\newcommand{\EG}{{\it e.g.}, }

\newcommand{\qed}
{\hspace*{\fill}$\square$\\}
\newcommand{\bul}{$\bullet\,$}
\newcommand{\ds}{\displaystyle}

\newcommand{\solS}{{\mathcal S}_{\lambda, \varphi}}
\newcommand{\solZ}{{\mathcal Z}_{\lambda, \varphi}}
\newcommand{\cbeta}{\check{\beta}}
\newcommand{\Ker}{{\mathcal K}_{\lambda, \varphi}}
\begin{document}
\begin{center}{
{\large {\bf
Reeb components with complex leaves 
\\
and their symmetries I : }
\vspace{7pt}\\
{
The automorphism groups 
\vspace{3pt} \\
and Schr\"oder's equation on the half line}
}
\vspace{12pt}
\\ 
Tomohiro HORIUCHI and Yoshihiko MITSUMATSU\footnote
{
The second author was partially supported by 
Grant-in-Aid for Scientific Research (B) No. 22340015. 
}
}
\end{center}
\begin{abstract}
We review the standard Hopf construction of Reeb components 
with leafwise complex structure 
and determine the group of 
leafwise holomorphic smooth automorphisms 
for tame Reeb components 
in the case of complex leaf dimension one. 
For this, we solve
the Schr\"oder type functional equation on the half line 
for expanding diffeomorphism.  
As a result, we see that the automorphism group of 
one with trivial linear holonomy on the boundary 
contains an infinite dimensional vector space, 
while in the case of non-trivial linear holonomy 
the group is of finite dimensional.  
\end{abstract}
\setcounter{section}{-1}
\section{\large Introduction}

The aim of this article is to begin a study of 
Reeb components in foliation with comlpex leaves of codimension one, 
especially focused on the symmetry in the real $3$-dimensional case. 

Quite often we call them leafwise complex foliations.

Recall that a $(p+1)$-dimensional Reeb component is 
a compact manifold $R=D^p\times S^1$ with a (smooth) foliation 
of codimension one, whose leaves are graphs of smooth functions 
$f : \mathrm{int} D^p \to \R$ 
where $\lim_{z \to \partial D^p}f(z) =+\infty$, 
and a compact leaf which is the boundary $S^{p-1}\times S^1$.   
Here we identify $R$ with $(D^p \times \R) / \Z$.  
See also the figures in Section 2.

Foliations of codimension one with complex leaves are 
drawing attentions in several complex variables because 
it appears as the Levi foliations of 
Levi-flat real hypersurfaces in complex manifolds. 
A simple construction of Hopf manifolds admits 
a Levi-flat real hypersurface, whose Levi foliation 
consists of a pair of Reeb components. 
This construction is generalized in Nemirovskii's examples \cite{Ne}. 
They have non-trivial linear holonomy along toral leaves. 
In this paper, from rather topological points of view, 
we study Reeb components with all kinds of holonomy.  
Of course, 
the case where the holonomy is {\it flat to the identity}, 
\IE it is infinitely tangent to the identity at the origin,  
is included. 
Such Reeb components appear in turbulization of a codimension one 
foliation along a closed transversal, or in pasting constructions 
like Dehn surgery of 3-manifolds.  


The present paper is organized as follows. 
In Part I we study Reeb components with 
complex leaves and their symmetries, 
as well as leafwise complex foliations in general. 
Part II is devoted to the study of 
Schr\"oder's functional equation on the half line, 
whose results are the core part of the study of 
the automorphism group of Reeb components. 

In Section 1 we fix the basic notions on foliations 
with complex leaves. They are also called 
{\it leafwise complex foliations}.   
After these preparations, Reeb components and the turbulization are 
reviewed in the context of leafwise complex foliations
in Section 2. 
Here we also review the notion of {\it tameness} of Reeb components, 
which we always assume in studying the symmetries.

Then, we study the symmetries of a 
Reeb component with complex leaves on a 3-manifold.

In section 3 we investigate the structure of 
the group of automorphisms of 
Reeb components with complex leaves of 
complex-dimension 1. 
Except for the case where the centralizer 
of the holonomy diffeomorphism in 
${\mathit Diff}^\infty [0, \infty)$ is 
not exactly known, 
we completely determine the structure of 
the automorphism group.  
This exception happens for some of 
diffeomorphisms which is flat to the identity. 
In anyway, we see that the automorphism group 
is of finite dimensional if the linear holonomy 
of the boundary leaf is non-trivial.  
On the other hand, if the linear holonomy is 
trivial, the automorphism group always 
contains an infinite dimensional vector space.  
(Theorem \ref{MainThm}).  
Such a clear contrast results from the 
analysis of Schr\"oder's equation 
in Part II. 
Similar results are obtained for Reeb components of 
complex leaf dimension $2$.  
They are explained by one of the authors in \cite{Ho}.

In Section 4 some direct corollaries to the 
results in Section 3 are stated. 
For example, the automorphism group of 
a Reeb foliation with complex leaves on the 
three sphere is understood.

The study of moduli space of tame Reeb component 
is studied by Meersseman and Verjovsky \cite{MV}.   
The moduli exhibits to a certain degree a similar phenomena 
to those of compact complex manifolds, 
especially concerning the finite dimensionality.  .  
As to automorphism groups 
our result tells that 
only the Reeb components with non-trivial 
linear holonomy on the boundary shows such a similarity.


The second part of the paper is devoted to 
the study of 
Schr\"oder's equation on the half line.   
It is in a form which is looking for 
eigen solutions for a pull-back operator. 
Here the pull-back diffeomorphism is nothing but 
the holonomy of the Reeb component when it is 
applied to Part I. 
In fact the results are the main ingredients 
in describing the automorphism groups 
of Reeb components in Part I. 
In section 5,  
we describe the space of solutions 
to 
Schr\"oder's equation.  
We also extend the values of the equation 
to $\C^2$ or still higher dimensional case, 
which is used 
in \cite{Ho}.  

For a diffeomorphism of the half line with 
non-trivial linear part, the computation is 
easy and in fact is well-known.  
For expanding diffeomorphisms with trivial 
linear part at the origin,  
the proofs are given in the subsequent sections. 
We present two different proofs.  

In Section 6 
a direct proof is given 
for the diffeomorphisms which are flat to the identity.  
This proof has a similar flavor to 
one by the center manifold theory, 
which is given in the final section. 

In Section 7, 
a proof given 
for 
diffeomorphisms with non-trivial finite jets, 
\IE those with the Taylor expansion at the origin 
different from the identity. 
This proof relies on 
Takens' normal form \cite{Ta} and the classical 
Fourier series on smooth function on the circle. 

Neither of the proofs in Section 7 nor 8 
works in other cases.

Then in the final section, 
we give a unified proof 
which is applicable for both of the above cases.   
The main tool is the center manifold theory 
of partially hyperbolic dynamical systems.  
This proof was suggested by Masayuki Asaoka. 

Throughout this article, we assume manifolds and 
foliations to be smooth unless otherwise stated.

The authors are deeply grateful to the members of Saturday Seminar 
at TITech, especially to Takashi Inaba, 
for exciting discussions and valuable comments, 
as well as to Masayuki Asaoka for suggesting us 
a proof in Part II and 
for introducing us the background. 
They also would like to express their gratitude 
to Laurent Meersseman and Alberto Verjovsky 
for their guidance on the basic materials in 
Section 1 and 2. 
\vspace{30pt}

\noindent
{\bf \Large  Part I : 
Reeb components with complex leaves}
\section{\large Basic definitions} 
Let $M$ be a $(2n+q)$-dimensional smooth manifold 
and $\F$ be a smooth foliation of codimension $q$ 
on $M$ and let $p=2n$ be the dimension of leaves. 
In this section and in the next, 
$n$ and $q$ are reserved 
for the complex leaf dimension and the codimension. 
We refer general basics for foliation theory to \cite{CC}. 

\begin{df}[Leafwise complex foliation, cf. \cite{MV}]
\label{LeafwiseComplexStructure}
{\rm \quad A smooth foliation $\F$ on a smooth manifold 
$M$ is said to be a 
{\it leafwise complex foliation} 
or 
{\it foliation with complex leaves} 
if there exists a system of local smooth foliated coordinate charts 
$(U_\lambda, \varphi_\lambda)$ 
where $\varphi_\lambda : U_\lambda \to V_\lambda \subset 
\C^n\times \R^q=\{(z_1,\, \cdots, z_n, y_1, \,\cdots, y_q)\}$
is a smooth diffeomorphism onto an open set $V_\lambda$ 
such that the coordinate change 
$(w_1,\, \cdots, w_n, t_1, \,\cdots, t_q)=
\gamma_{\mu \lambda}(z_1, \,\cdots, z_n, y_1, \,\cdots, y_q)$ 
is smooth, $t_j$'s depend only on $y_k$'s ($j,k =1, \,\cdots , q$), 
and when $y_k$'s are fixed $w_l$'s are holomorphic in $z_m$'s, where 
$\gamma_{\mu \lambda}: \varphi_\lambda(V_\lambda \cap U_\mu)
\to \varphi_\mu(V_\lambda \cap U_\mu)$. 
}
\end{df}

\begin{rem}{\rm \quad 
Instead of assuming local coordinate 
system as above, it is also natural to 
consider the following condition 
that 
the smooth foliation $\F$ admits 
a smooth almost complex structure $J$ acting on 
the tangent bundle $T\F$ to the foliation, 
which is integrable on each leaves, 
namely there exists local holomorphic coordinates on each leaves.  
We assume here $J$ is smooth on the ambient manifold $M$. 
(Of course $J$ becomes more than smooth in each leaf. )

This might appear slightly weaker than Definition 1.1,  
but 
eventually they are equivalent to each other. 
To prove from the weaker to the stronger 
is nothing but 
the parametric version of Newlander-Nirenberg's theorem. 
The Newlander-Nirenberg theorem \cite{NN} in the usual sense 
claims that an almost complex manifold $(L,J)$ 
admits a complex structure 
(a holomorphic local coordinate system) if 
the Nijenhuis tensor $N_J$ vanishes. 
In dead in \cite{NN} 
Newlander and Nirenberg mentioned 
in the very last paragraph 
that the parametric version holds. 
For the case of $n=1$, 
even the parametric version seems to be classically known, 
 \EG see \cite{Mo}.

It is also well known that if an almost complex structure 
$J$ is real analytic 
on $2n$-dimensional real analytic manifold $L$, 
the Newlander-Nirenberg theorem has a simple geometric proof. 
See, for example, Appendix A4 of \cite{Hu}.  
Using this argument, 
if we can take such a smooth foliated chart 
$(u_1, u_2, \cdots, u_{2n}, y_1, \cdots , y_q)$ 
that $J$ is real analytic on $(u_1, \cdots, u_{2n})$, 
we can show the existence of a local coordinate system 
in Definition 1.1. 

It should be also remarked that 
the case of Levi-flat real hypersurface $M$ in 
an $(n+1)$-dimensional complex manifold $W$, 
the stronger one is easily satisfied. 
(If $M$ is of class $C^r$, then we can only assure that 
$TM$ is of class $C^{r-1}$, 
so that the resultant local coordinate system is assured 
to be of class $C^{r-1}$. )  
}
\end{rem}

\begin{df}{\rm \quad 
A diffeomorphism between two foliated manifolds 
with complex leaves is said to be an {\em isomorphism} 
between leafwise complex foliations iff it preserves the foliations 
and gives rise to biholomorphisms between leaves. 
An {\it automorphism} is an isomorphisms between the same one. 
}
\end{df}

In this paper we are mainly concerned 
with foliation of codimension one. 
In particular, 
our interest will be focused on Reeb components of real dimension $3$, 
namely in the case of $n=1$ and $q=1$. 
As we see from the examples of Nemirovskii \cite{Ne} 
even a real analytic Levi-flat hypersurface in a complex manifold 
can admit Reeb components in its  Levi foliation.  
In such a case, the holonomy along the toral boundary leaf 
has a non-trivial linear part. 

Apart from Levi-flat real hypersurfaces, 
for example, if we perform a turbulization 
we easily find various leafwise complex foliations  
admitting Reeb components with holonomy 
flat to the identity. 
See the next section for more detail.

\section{\large Reeb components with complex leaves}
In this section we review a particular construction of 
Reeb component with complex leaves and a process of turbulization 
which produces a new Reeb component in a leafwise complex foliation. 


In order to make pasting construction easier, 
we introduce the following notions. 
Let $(R, \F, J)$  (or simply $R$ for short) 
be a Reeb component with leafwise complex structure 
of complex leaf dimension $n$ and $(H,J_H)$ be its boundary 
leaf. 
\begin{df}\label{TameFlat} 
{\rm 
The Reeb component $R$ has a {\it tame} boundary 
(or `$R$ is {\it tame} at boundary' for short, 
or even shorter `$R$ is {\it tame}') 
with respect to a product coordinate $H\times[0,\varepsilon)$ 
of a collar neighborhood of $H$ 
if it gives rise to a smooth foliation with 
leafwise complex structure 
when it is  pasted with the product foliation 
$(H\times(-\varepsilon, 0], 
\{H\times\{x\}\vert x\in(-\varepsilon, 0]\}, J_H)$ 
along their boundary. 
Here each leaf $H\times\{x\}$ has the same complex structure 
as $H$ when identified with the natural projection. 
Namely, 
the Reeb component is extended to the outside as a product foliation. 
}
\end{df}
The notion of tameness was introduced in \cite{MV}.  

\begin{rem}{\rm \quad 
If we forget the leafwise complex structure and consider 
the same notion only as foliation of codimension one,  
it does not depend on the choice of product coordinate 
on the positive side and the tameness implies 
exactly that the holonomy is tangent to the identity 
to the infinite order. This is because the set of expanding 
diffeomorphisms of the half line $[0,\infty)$ 
which are infinitely tangent to the identity is 
an open convex cone and invariant under conjugation 
by any diffeomorphism.  
Also remark that the tameness depends only 
on the smooth projection of the collar neighborhood 
to the boundary, which the product coordinate defines. 
If two projections have the same infinite jets on the boundary, 
the tameness notion coincides for the two. 
}
\end{rem}

\begin{df}\label{simple}
{\rm 
The leafwise complex structure of a Reeb component $R$ is {\it simple} 
around boundary 
(or $R$ has a {\it simple} complex structures around boundary) 
if the boundary has a collar neighborhood $U\cong H\times [0, \varepsilon)$ 
such that 
the restriction of the projection $U=H\times [0, \varepsilon) \to H$ 
to each leaf in $U$ is holomorphic.  

This notion should also be understood relative to the projection 
from a collar neighborhood to the boundary.  }
\end{df}

The notions of tameness and simpleness apply not only to Reeb components 
but also to more general leafwise complex foliations 
of codimension one with a compact leaf or a boundary leaf.  

Clearly if a Reeb component has simple complex structures around the boundary 
and the holonomy of the boundary leaf is infinitely tangent to the identity, 
it is tame with respect to the appropriate projection. 
The tameness condition prohibits 
unexpected wild behaviour around boundary.  
In particular in the case of complex leaf dimension $=1$, 
it induces a strong consequence due to Meersseman and Verjovsky.  
 See the following subsection. 

\subsection{Reeb component by Hopf construction}
Let us recall the Hopf construction 
which is one of the standard ways to construct 
Reeb components. 
This construction gives rise to a tame Reeb component 
if the holonomy $\varphi$ is infinitely tangent to 
the identity at $\,x=0\,$.

\begin{constr}\label{Hopf-1}{\rm (Hopf construction)\quad 
Let $\varphi\in\mathit{Diff}^\infty([0,\infty))$ 
a diffeomorphism of the half line $\R_{\geq 0}=[0, +\infty)$ 
satisfying 
$\varphi(x)-x>0$ for $x>0$, 
namely the origin is an expanding unique fixed point.  
Also take a (local) biholomorphic diffeomorphism 
$G\in \mathit{Diff}^\mathit{hol}(\C^n, O)$ 
which is expanding.  This implies that for some 
small neighborhood $D$ of the origin $O$ with smooth boundary 
$G(\mathrm{int}D)\supset \overline{D}$ and 
$\ds \cap_{k=1}^{\infty}G^{-k}(D)=\{O\}$. 
Now take $\ds U= \cup_{k=1}^\infty G^k(D) \subset \C^n$.

Then on $\tilde{R}=U\times[0,\infty)\setminus
\{(O,0)\}\subset \C^n \times \R$, take the restriction $\tilde\F$ of the 
product foliation $\{\C^n\times\{x\}\}$ 
together with the natural complex structure on leaves 
and a diffeomorphism 
$T=G\times\varphi$ on $U\times[0,\infty)\setminus
\{(O,0)\}$). Practically we take fairy simple diffeomorphisms 
such as linear maps as $G$ 
so that $U$ becomes the whole $\C^n$. 
Then on the quotient $R=\tilde{R}/T^\Z $ a foliation $\F$ with 
complex leaves is naturally induced. 

From the construction, it is simple around the boundary.  
If the holonomy is infinitely tangent to the identity 
it is also tame with respect to the coordinate in the construction. 

The boundary $U\setminus\{O\}/G^\Z$ is a complex manifold which is 
a so called {\it Hopf manifold}. 
In the case $n=1$ it is an elliptic curve and the construction is 
equivalent to one with linear map as $G$. 
}
\end{constr}

\begin{rem}{\quad \rm It is well known as an elementary fact 
in complex dynamical systems 
that assuming $G(\mathrm{int}D)\supset \overline{D}$ 
for a bounded connected domain $D\subset \C^n$ is enough to conclude that 
there is a unique linearly expanding fixed point 
in $D$ and $D$ is included in the attracting basin 
of $G^{-1}$.  
}
\end{rem}

\begin{thm}\label{MV}
{\rm
(Meersseman-Verjovsky, \cite{MV})\quad 
Any tame Reeb component with complex leaves 
of complex dimension $1$ 
is isomorphic to 
one of those given by the Hopf construction.  
}
\end{thm}

We present a couple of extensions 
(variants) of the above 
construction. 

\begin{constr}\quad{\rm 
Now, let us take the product not with 
the half line but with the whole real line $\R$. 
Let $M$ and $\Phi\in \mathit{Diff}^\infty_+
(\R)$ be as follows. 
\vspace{5pt}
\\
\qquad
$\bullet$ $ M=(U\times\R\setminus\{(O,0)\})/T'^\Z, \quad
T'=G\times\Phi$, 
\vspace{3pt}
\\
\qquad 
$\bullet$  $x=0$ is an expanding unique fixed point of $\Phi$.  
\vspace{5pt}
\\
$M$ consists of two Reeb components 
and in exactly the same way as above    
a foliation with leafwise complex structure is induced on $M$. 

Note that in this construction and in the previous one, 
the holonomy of the toral leaf is given by $\Phi$ 
and by $\varphi$ respectively. 
}
\end{constr}

\begin{constr}\label{HopfSurface}\quad{\rm 
Next, we further extend the previous consruction to 
obtain a Reeb component as a part of a Levi-flat 
hypersurface in a Hopf surface. 
We take 
$\Phi$ to be a 
linear expansion in order to extend 
it a biholomorphic (in fact linear) expansion 
$\tilde\Phi$ of $\C$.   
Note that $\R$ is an invariant subspace in $\C$. 
Fix the expansion ratio $\mu >1$ 
and take the followings; 
\begin{eqnarray*}
W=(U\times\C\setminus\{O\})/T''^\Z, \,\, \quad\quad 
T''=G\times\tilde\Phi, \quad 
\tilde\Phi(z)=\mu z  \,\,\,(z\in\C), 
\\
M=(U\times\R\setminus\{(O,0)\})/T'^\Z, \quad 
T'=G\times\Phi, \quad \, 
\Phi= \tilde\Phi\vert_\R\, .\,\, \qquad \qquad
\end{eqnarray*}
$W$  is an $(n+1)$-dimensional Hopf manifold, 
$M$ is its Levi-flat real hypersurface 
with the Levi foliation 
consisting of two Reeb components, and 
a unique compact leaf is the Hopf manifold $(U\setminus\{O\})/G^\Z$ 
of $\dim_\C =n$.  
}
\end{constr}

\begin{prob}\quad{\rm 
Theorem \ref{MV}
due to Meersseman and Verjovsky 
poses the following questions. 
We assume the complex leaf dimension to be one.  
Provided that two Reeb components 
with leafwise complex structures 
have the same boundary holonomy and their boundary leaves 
are biholomorphic to each other, are they isomorphic 
as leafwise complex foliations? 
Does there exist a Reeb component with complex leaves 
which is not isomorphic to a tame one but with 
holonomy infinitely tangent to the identity? 
Or does there exist one which is not 
isomorphic to any of those given 
by the Hopf construction? 
One more similar but subtle question is to ask whether if 
a tame Reeb component is always isomorphic to 
one given by the Hopf construction.  

The second form of question seems less difficult and 
negative.  
Anyway, those questions are asking what should be the 
complete invariants to determine Reeb components 
without assuming the tameness.  }
\end{prob}

\begin{constr}\label{preparation}\quad
{\rm We introduce one more 
construction, 
which is a preparation for turbulization. 
Take $\tilde{M}= 
(\C^n\times \R) \setminus\{O\}\times(-\infty, 0]$ and restrict 
the product action $\hat{T}=G \times\Psi$ to $\tilde{M}$,   
where $\Psi$ is an orientation preserving diffeomorphism 
of $\R$ which fixes $0$,  
expanding on $[0,\infty)$, 
and {\it contracting} on the negative side 
$(-\infty, 0]$, 
\IE 
$\Psi(x) > x $ for $x<0$. 
On $\tilde M$ we take (the restriction of) 
the horizontal foliation $\tilde \F$.  
Then take the quotient 
$(M,\F, J_\F) = 
(\tilde M,\tilde \F, J_\mathrm{std})/\hat T ^\Z$. 

The non-negative part is nothing but the Reeb component 
constructed in \ref{Hopf-1} 
regarding $\varphi=\Psi\vert_{[0,\infty)}$. 
The non-positive side $(N,\G)=(M,\F)\vert_{x\leq 0}$ 
remains non-compact and 
is in fact a foliated 
$(-\infty, 0]$-bundle with holonomy 
$\psi= \Psi\vert_{(-\infty, 0]}$. 

If we remove the boundary compact leaf $\{x=0\}$ 
from the non-positive side 
$(N, \G)$, 
it is isomorphic to $(\C^n\setminus \{O\})\times S^1$. 
For a better description of turbulization process, 
let us be more precise about this identification. 
This is done by embedding $\hat T$ in a 1-parameter family. 
Take a smooth curve $G_t$ in $\mathit{GL}(2;\C)$ 
and also a smooth curve $\psi_t$ in 
$\mathit{Diff}^\infty((-\infty, 0])$ 
satisfying the following conditions. 
$$
\psi_k=\psi^k \,\,\, (k\in \Z), \quad 
\psi_{t+1}=\psi\circ\psi_t \,\,\, (t\in \R), \quad 
\frac{\partial\psi_t(x)}{\partial t}>0\,\, (\forall x, t),
$$
$$
G_k=G^k \,\,\, (k\in \Z), \quad 
G_{t+1}=G\circ G_t \,\,\, (t\in \R). 
$$
Then, fixing (any) $x_0<0$,   
$x=\psi_t\,(x_0)$ gives a diffeomorphism between 
$(-\infty, 0)\,(\ni x)$ and $\R(\ni t)$.   
Then the identification of 
$(z,x)\in(\C^n\setminus\{O\})\times (-\infty, 0)$ with
 $(w,t)\in(\C^n\setminus\{O\})\times \R$  
by 
$(z=G_t(w), x=\psi_t(x_0))$ 
conjugates $\hat T\vert_{x<0}$ into $(w,t)\mapsto(w, t+1)$. 

Of course an easy way to choose such a 1-parameter family 
$\,\{\psi_t\}\,$ is to take a 1-parameter subgroup. 
Take a smooth vector field $\,\rho(x)\frac{d}{dx}\,$ on 
$\,(-\infty, 0]\,$ with $\,\rho(x)>0\,$ for $\,x<0\,$ and 
$\,\rho(0)=0\,$. 
Then putting $\,\psi_t=\exp(tX)\,$ we obtain such 
$\,\psi_t\,$ with $\,\psi=\psi_1\,$. 
If we choose $\,\rho(x)\,$ to be flat at $\,x=0\,$, 
$\,\psi_t\,$ is infinitely tangent to the identity at $x=0$. 
Even if $\,\psi\,$ does not ly 

It is worth remarking that this identification gives rise to 
a partial compactification of horizontally foliated manifold 
$((\C^n\setminus\{O\})\times S^1, 
\{(\C^n\setminus\{O\})\times\{t\}\}$ 
by a Hopf manifold $N$ as a boundary leaf so as to 
obtain $(N, \G)$. 
If we take diffeormorphisms $\psi_t$ 
infinitely tangent to the identity at the origin, 
we obtain a tame structure. 
Also on the non-negative side, 
by taking smilar family $\varphi_t$ 
for $\varphi=\varphi_1$, we also obtain a tame structure on 
the non-negative side.  
Once we obtain tame ones on both side with the same 
complex structure on the boundaries, we can paste them 
to obtain a smooth structure.  
}
\end{constr}

\subsection{Turbulization 
in $L\times S^1$}\label{standard turbulization}
Here we review the turbulization, 
which is classically well-known modification of 
a foliation of codimension one 
to yield a new Reeb component.  
We start from a standard situation. 

\begin{constr}\label{standard turbulization}{\rm \quad
Let $(M,\F)$ be a leafwise complex foliation of codimension one 
and assume that there is an embedded 
solid torus $U={\mathrm{int }}D^{2n} \times S^1$ 
on which the the induced foliation 
is $\{{\mathrm{int }}D^{2n} \times \{*\}\}$ and 
the induced complex structure is also the canonical ones 
on each 
${\mathrm{int }}D^{2n} \times \{*\}
\cong {\mathrm{int }}D^{2n}$ 
$\subset \C^n$.  Let $(w,t)$ 
denote the natural coordinate of 
$U={\mathrm{int }}D^{2n} \times S^1$ 
where $S^1$ is regarded as $\R/\Z$. 
Then we remove 
$\{O\}\times S^1$ from $U$ 
and let 
$U^*$ 
denote the result. 
Using the coordinate 
$(w,t)$ $U^*$ 
is identified with an open subset 
of the negative side of Construction \ref{preparation}, 
together with leafwise complex foliations.  
Therefore we can compactify this end 
with the Hopf manifold as in Construction \ref{preparation} 
and also 
if we add positive side of Construction \ref{preparation} 
we obtain a leafwise complex foliated manifold 
without boundary 
but with a new Reeb component. 
For this construction 
we can choose any of $G_t$, $\psi_t$, and $\Phi$ as in 
Construction \ref{preparation}.    
The above process including adding the positive side is 
the leafwise complex version of the {\it turbulization}.  
See also Figure 1 below.  
}
\end{constr}

\subsection{General case}
It is easy to find a closed transversal to a foliation of codimension one, 
namely, an embedded circle which is transverse to the foliation, 
unless the manifold is open and the foliation is too simple.  
Like in the case of 
a smooth foliation without leafwise complex structure, 
it is always possible to perform the turbulization 
in a tubular neighborhood of any closed transversal 
regarding leafwise complex structure.  
This fact also belongs a kind of folklore, while below 
it is reviewed.  
\begin{thm}{\rm
\quad 
Let $(M^{2n+1}, \F, J)$ be a smooth leafwise complex foliation 
of codimension one and $K\subset M$ is a closed transvesal, 
namely there exists a smooth embedding $f:S^1 \to M$ 
which is transverse to the foliation $\F$ with 
its image $f(S^1)=K$. 

Then, there exists a tubular neighborhood 
$U\cong K \times {\mathrm{int}\,}D^{2n}$ such that 
the restricted foliation $(U, \F_U, J\vert_{\F_U})$ 
is isomorphic to the standard one 
$(S^1 \times {\mathrm{int}\,}D^{2n}, 
\F_0=\{t\}\times {\mathrm{int}\,}D^{2n}, J_0)$ 
and through this isomorphism 
$K$ is identified with $S^1\times \{O\}$. 

In particular, we can perform the standard turbulization 
\ref{standard turbulization} 
in $U$. 
}
\end{thm}
This theorem is a direct corollary to the following lemma.  
\begin{lem}
{\rm
\quad
The group 
$\mathit{Diff}^{\mathit{hol}}(\C^n, O)$
of germs of holomorphic diffeomorphisms 
of $(\C^n,  O)$ which fix the origin is pathwise connected.  
}
\end{lem}
The lemma immediately follows from the two facts 
that $GL(n;\C)$ is pathwise connected 
and that such a germ with identical linear part can 
be joined by a straight segment to the identity.

\subsection{Dehn surgery in $\dim=3$ vs. higher dimensional turbulization}
In order to close the section, 
this subsection provides with some remarks concerning the possibility 
of pasting the Reeb component in a different way in a turbulization. 
In the rest of this section, 
we assume the holonomy $\Psi$,  and eventually 
$\varphi$ and $\psi$, 
to be infinitely tangent to the identity at the origin,  
so that it is easier to past two peices along their boundarie. 

\begin{rem}{\rm \quad 
If we forget the leafwise complex structure and treat foliations 
only as smooth objects, basically there are two ways to perform 
the turbulization. The one has been already described above and 
is indicated in Figure 1. 
For the other one we can reverse the top and bottom 
of the Reeb component (Figure 2). This is because the cyclic 
(universal for $n\geq 2$) covering of the boundary leaf is 
$\R^{2n}\setminus\{O\}\cong S^{2n}\setminus\{N, S\}$ and 
two ends are exchangeable by a diffeomorphism. 
However, as a complex manifold, $\C^n\setminus\{O\}$ has 
one convex end and the another concave.  
For the case of complex leaf dimension $n$ greater than one,  
these two ends are not exchangeable.  
In particular, for $n\geq 2$, 
the turbulization for leafwise complex foliations does not 
change the homotopy class of the tangent bundle.

}
\end{rem}
\noindent
\begin{figure}[H]
 \begin{center}
  \includegraphics[width=6cm]{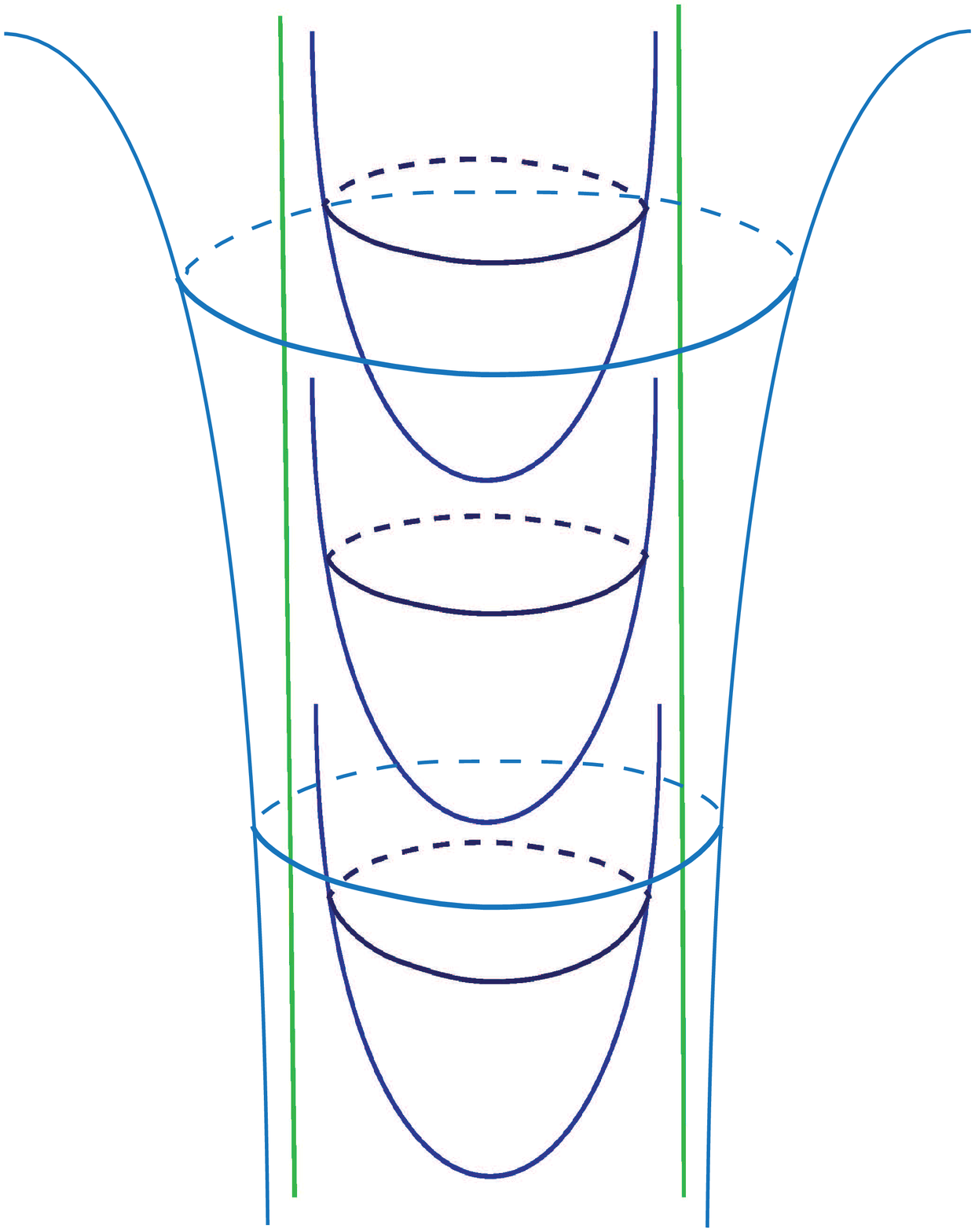}
\qquad
  \includegraphics[width=6cm]{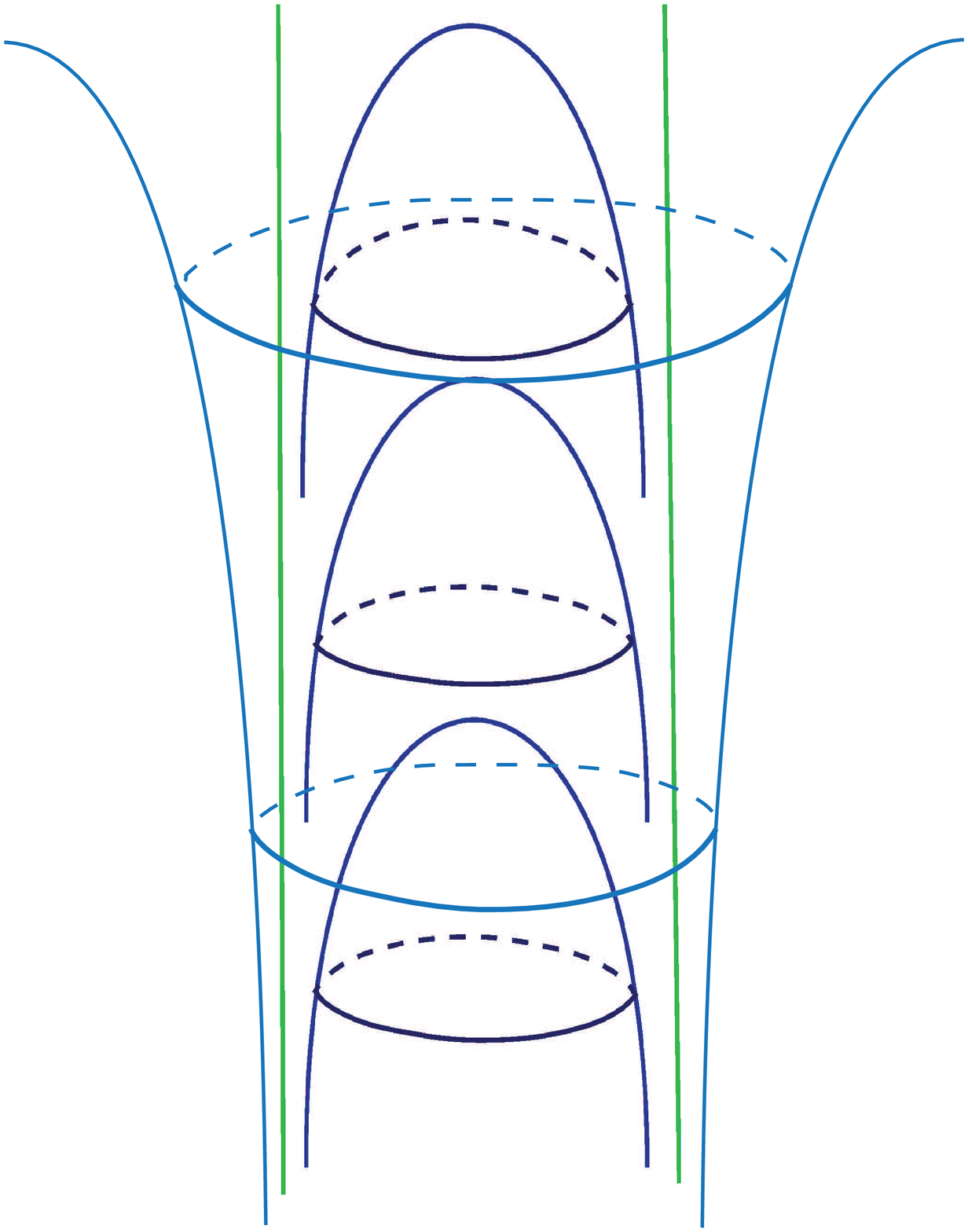}
\\
\hspace*{\fill}
Figure 1
\hspace*{\fill}
\hspace*{\fill}
Figure 2
\hspace*{\fill}
\\
\hspace*{\fill}
\begin{minipage}{119mm}
{\small \mbox{}
\\
The green lines indicate the bounndary 
leaves of Reeb components. 
The axes of the rotational 
symmetries of the Reeb components, 
which are not drawn in the figure, 
correspond to $\{O\}\times \R_+$. 
}
\end{minipage}\vspace{-10pt}
\hspace*{\fill}
 \end{center}
\end{figure}

\begin{rem}{\rm \quad In the case of complex leaf dimension one, 
the `upside-down' construction always works.  
Namely, in Construction \ref{Hopf-1}, $z \leftrightarrow z^{-1}$ 
always induces an biholomorphism on the boundary elliptic curve.  
Therefore we can regard it as a turbulization inthecated in Figure 2. 

If the boundary elliptic curve admits a complex multiplication, 
namely finite but discrete symmetries of order 2, 3 or 4, 
removing the Reeb component and pasting it back with one of 
those symmetries is a special kind of Dehn surgeries. 

More generally, in the process of turbulization, 
after removing a tubular neighborhood of 
the closed transversal and compactify the new boundary 
(namely after the process of Construction \ref{preparation}), 
instead of filling up with a Reeb component as explained above, 
we fill up the boundary as follows.  
We prepare another 
Reeb component with a different complex structure. 
If their boundaries match up through some diffeomorphism,    
we can fill up the boundary with that Reeb component. 
In this way, a Dehn surgery corresponding to 
any element of the mapping class group 
${\mathcal M}_1(\cong \mathit{}SL(2;\Z))$ 
of a $2$-dimensional torus $T^2$ is realized 
for a closed transversal in a leafwise complex 
codimension one foliation of $n=1$.  
}
\end{rem}

\section{\large Symmetries of $3$-dimensional Reeb components }
In this section we compute the group of automorphisms of 
a Reeb component of dimension $3$ with complex leaves, 
which is given by the Hopf construction.  
In order to fix notations, 
we present our objects again. 
Let $\tilde R$ be $\C\times[0,\infty) \setminus \{(0,0)\}$, 
take $\lambda\in \C$ with $\vert\lambda\vert>1$ and 
a diffeomorphism $\varphi\in\mathit{Diff}^\infty([0,\infty))$ 
which is expanding, namely, satisfying 
$\varphi(x)-x>0$ for $x>0$.   
Let $G$ denote the linear expansion of $\C$ defined 
as the multiplication by $\lambda$ 
and $T:\tilde{R} \to \tilde{R}$ be $T=G \times \varphi$.    
Then we obtain a Reeb component 
$(R, \F, J)= 
(\tilde{R}, \tilde{\F}, J_\mathrm{std}) /T^\Z$ 
as the quotient, as well as the boundary elliptic curve 
$H=\C\setminus\{0\}/G^\Z$.   
Here, on the upstairs the leaves of the foliation $\tilde{\F}=
\{\C\times\{x\} \vert x>0\}\sqcup \{\C^*
\times\{0\}\, \}$ 
are equipped with the natural complex structure 
$J_\mathrm{std}$ which is inherited by those of $\F$. 
\vspace{8pt}

We classify the diffeomorphisms $\varphi$ 
into the following three cases 
according to the nature of its jet at $x=0$.  
Here, $\varphi^{(i)}$ denotes the $i$-th derivative of 
$\varphi(x)$ 
and $j^i\varphi(0)$ denotes the $i$-th jet at $x=0$. 
\vspace{5pt}
\\
\noindent
\quad Case (1) : $\varphi'(0)=\mu >1$ , 
\vspace{2pt}
\\
\noindent
\quad Case (2) : for some $n \geq 2$ ,   
$j^{n-1}\varphi(0)=j^{n-1}\mathit{id}(0)\,$  
and 
$\,\,\varphi^{(n)}(0)>0$ ,  
\vspace{2pt}
\\
\noindent
\quad Case (3) : 
$j^\infty\varphi(0)=j^\infty\mathit{id}(0)\,$ .  
\vspace{8pt}

The discussions in Subsection \ref{LiftRest} 
does not depend on the above classification.  
As reviewed in Section \ref{centralizer}, 
the structure of the centralizer of $\varphi$ 
is very subtle for a certain class in Case (3). 
For the rest of Case (3) and for Case (1) and (2), 
the centralizer is fairy simple.  

The main result of this paper is the 
computation of the automorphisms 
which fixes the boundary and 
the transverse space.  
This is in deed the main results of Part II 
of this paper.  
Concerning this part, for Case (2) and for Case (3) 
the results are the same. 
For Case (1), such automorpsims are very few. 

The internal structure 
of the automophism group 
for Case (2) and (3) rather depends only on 
the nature of the centralizer 
(see Subsection \ref{structure}).   
The centralizer itself is included in the 
automorphism group. 
Except for this part,  
the extendability to the outside is basically the same 
for any $\,\varphi\,$ in Case (2) and (3).

\subsection{Lift to $\tilde R$ and restriction to $H$}
\label{LiftRest}
Let us consider the group $\mathit{Aut} (R, \F, J)$, which 
is also denoted by $\mathit{Aut} R$ for short, of 
all foliation preserving diffeomorphisms of $R$ whose 
restriction to each leaf is holomorphic.  
Also we consider the group of holomorphic diffeomorphisms 
$\mathit{Aut} H$ of the boundary elliptic curve $H$ 
as well as its identity component $\mathit{Aut}_0 H$ which 
is isomorphic to $T^2$ and can be identified with $H$ itself.

\begin{prop}{\rm
\quad The image of the restriction map 
$r_H:\mathit{Aut} R \to \mathit{Aut} H$ is 
exactly  $\mathit{Aut}_0 H$.  
}
\end{prop} 
{\it Proof.} \quad 
If we regard 
$\mathit{Aut} H / \mathit{Aut}_0 H$ 
as a subgroup of $\mathit{SL}(2;\Z)$, 
in most cases it is just $\{\pm E\}$ where $E$ denotes 
the identity matrix. 
In a few cases where the elliptic curve $H$ admits 
complex multiplications, they are of order 3,4, or of 6 
and a kind of `rotations' on the universal covering, 
\IE elliptic matrices in $\mathit{SL}(2;\Z)$. 
In any of those cases, no element in 
$\mathit{Aut} H \setminus \mathit{Aut}_0 H$ 
preserves the direction of holonomy and thus none 
extends to $R$ as a foliation preserving diffeomorphism.

On the other hand, any element in  $\mathit{Aut}_0 H$ 
is obtained as the quotient of the scalar multiplication 
$m_a : \C^* \to \C^*$ by some nonzero complex number $a$. 
The automorphism $m_a\times \mathrm{id}_{[0,\infty)}$ of  
$\tilde R$ clearly descends to $R$ and defines an element 
in $\mathit{Aut} R$. 
\hspace*{\fill}$\square$
\\
\\
By this proposition, the study of the structure of 
$\mathit{Aut} R$ breaks into two parts,  
that of the kernel $\mathit{Aut}(R, H)$ 
and the study of the restriction map $r_H$. 
\\

Now it is easier to look at the lifts of automorphisms 
on $\tilde R$. 
Any element $f \in \mathit{Aut} R$ has a 
lift $\tilde f\in\mathit{Aut}(\tilde R, \tilde \F, J_\mathrm{std})$ 
($=\mathit{Aut}\tilde R$) 
which takes the form 
$$
\tilde f(z,x)=(\xi(z,x), \eta(x))
$$
in $\C\times [0,\infty)$-coordinate. 
A lift $\tilde f$ should commutes with the covering 
transformation $T$, because, $T\circ \tilde f = \tilde f \circ T^k$ 
for some $k\in \Z$ but it is easy to see that $k=1$ 
when it is restricted to the boundary.  
Therefore an element in 
$\mathit{Aut}\tilde R$ 
is a lift of some element in $\mathit{Aut} R$ 
if and only if it commutes with $T$.  
Let 
$\mathit{Aut}(\tilde R ; T)$ 
denote the centralizer of $T$ in 
$\mathit{Aut}\tilde R$,  
namely, the group of all such lifts. 
It contains an abelian subgroup 
$\{m_a\times \mathrm{id}_{[0,\infty)}\vert a\in\C^*\}\cong\C^*$.   
This subgroup injectively descends to a subgroup of 
$\mathit{Aut}R$ 
which restricts exactly to 
$\mathit{Aut}_0 H \cong \C^*/\lambda^\Z$.  
It is important to remark that 
whether $\mathit{Aut}_0 H$ admits a homomorphic section 
is not a trivial question.  
Postponing this question until the end of this section, 
we go on an easier way.

Let us introduce one more subgroup 
$\mathit{Aut}(\tilde R, \tilde H; T)$ of 
$\mathit{Aut}(\tilde R ; T)$ 
which consists of all elements 
which act trivially on the boundary $\tilde H$.  
Any element 
$f \in \mathit{Aut}(R, H)$ has a unique lift to 
an element $\tilde f\in\mathit{Aut}(\tilde R, \tilde H; T)$ 
Namely,    
\begin{cor}{\rm \quad 
$\mathit{Aut}(R, H)$ 
is isomorphic to  
$\mathit{Aut}(\tilde R, \tilde H; T)$. 
}
\end{cor}
Again, let 
$\,g\in \mathit{Aut}\tilde R\,$ 
be presented 
in the form 
$\,
g(z,x)=(\xi(z,x), \eta(x))
\,$.  
\begin{lem}\label{principal lemma}
{\rm \quad 
The element $\,g\,$ 
in $\,\mathit{Aut}\tilde R\,$  
belongs to 
$\mathit{Aut}(\tilde R ; T)$ 
if and only if 
it satisfies the following conditions. 
\begin{itemize}
\item[(a)] $\xi(z,x)= a z + b(x)$, $b(0)=0$ 
for some 
$b\in C^\infty([0,\infty), \C)$ and $a\in\C^*$.    
\item[(b)] $b(\varphi(x))=\lambda b(x)$. 
\item[(c)] 
$\varphi\circ \eta = \eta \circ \varphi$, 
namely, $\eta \in Z_\varphi=$ the centralizer of 
$\varphi$ in $\mathit{Diff}^\infty([0,\infty))$.  
\end{itemize}
Further more, 
$g$ belongs to $\mathit{Aut}(\tilde R, \tilde H ; T)$ 
if and only if the above conditions are satisfied 
with $a=1$. 
}
\end{lem}
{\it Proof.}\quad 
Let us first show the {\it only if} direction,  
then the {\it if} direction will become almost trivial.  

Assume 
$g\in\mathit{Aut}(\tilde R ; T)$.   
$\xi(z,x)$ is smooth and holomorphic in $z$. If $x$ is fixed, 
$\xi(\,\cdot\,, x): \C
\to\C
$ is 
a holomorphic automorphism  
even in the case where 
$x=0$ 
because the origin is a removable singularity, 
it is a linear map with nontrivial linear term. 
Therefore it is written in the following form; 
$\xi(z, x)=a(x)z + b(x)$ where $a(x), b(x) 
\in C^\infty([0,\infty), \C)$ 
with $a(x)\ne0$ and $b(0)=0$.  
These also apply to elements in 
$\mathit{Aut}\tilde R$.  

Now look at the commutation relation 
$\,g\circ T=T\circ g\,$.  
This implies 
$$
(a(\varphi(x))\lambda z + b(\varphi(x)), \eta(\varphi(x))
=(\lambda a(x)z + \lambda b(x), \varphi(\eta(x)))\, .
$$
Thus we obtain (b) and (c).  
This also tells us that $a(\varphi(x))=a(x)$,  
so that for any $x\geq0\,$ we have   
$\ds a(x)=\lim_{n\to\infty}a(\varphi^{-n}(x))=a(0)$ 
and (a) is concluded. 

For $g\in\mathit{Aut}(\tilde R, \tilde H ; T)$ 
we just need to confirm that $a=1$. 
\hspace*{\fill} $\square$
\begin{rem}{\rm\quad 
The condition (b) 
appears as Equation (I) in Part II. 
Solving this 
Schr\"oder type functional equation 
on the half line $\,[0, \infty)\,$ for 
given $\,\lambda\,$ and $\,\varphi\,$ 
is the main theme in Part II. 
}
\end{rem}
\begin{cor}\label{identification}
{\rm \quad 
$\mathit{Aut} R$ is naturally isomorphic to 
$\mathit{Aut}(\tilde R; T)/T^\Z$. 
}
\end{cor}

\subsection{Centralizer of 
$\varphi$ in $\mathit{Diff}^\infty([0,\infty))$}
\label{centralizer}

For an expanding diffeomorphism 
$\,\varphi\in \mathit{Diff}^\infty([0,\infty))\,$,  
concerning its centralizer in 
$\,\mathit{Diff}^\infty([0,\infty))\,$
and the embeddability 
in a (smooth) 1-parameter subgroup, 
the followings are known. 
\begin{thm}{\rm\quad 
1) \quad For Case (1) thanks to Sternberg's 
linearization \cite{St} 
and 
for Case (2) thank to Takens' normal form \cite{Ta},  
there exists a smooth vector field 
$\ds\,X=\rho(x)\frac{d}{dx}\,$ on $\,[0,\infty)\,$ 
such that $\varphi=\exp{X}$ and the centralilzer 
$Z_\varphi$ exactly coincides 
with the 1-parameter subgroup $\,\exp(tX)\,; \, 
t\in \R\}\,$ generated by $X$.  
\\
2)\quad In Case (3), if there exists a smooth 
vector field 
$\ds\,X=\rho(x)\frac{d}{dx}\,$ on $\,[0,\infty)\,$ 
with $\varphi=\exp{X}$, 
then its centralizer  
$Z_\varphi$ in coincides 
with the 1-parameter subgroup 
$\,\exp(tX)\,; \, t\in \R\}\,$.  
}
\end{thm}

In general, the centralizer $Z_\varphi$ of 
$\varphi$ which is infinitely tangent to the identity 
at $\,x=0\,$ is known 
to be fairy wild (see \cite{Ey}).  
For an expanding  $\,\varphi\,$, 
  it is known that 
there exists a unique $C^1$-vector field 
$\ds X_\varphi=\rho(x)\frac{d}{dx}$ 
on the half line $[0,\infty)$, 
which is of class $C^\infty$ on $(0, +\infty)$, 
in such a way that 
the exponential map 
$\exp X_\varphi$ 
coincides with $\varphi$  (see \cite{Sz} and \cite{Na}).  
This vector field is often called 
the {\it Szekerez vector field} of $\varphi$. 
If the Szekerez vector field $X_\varphi$ is of class 
$C^\infty$ on $[0,\infty)$, 
namely $\rho(x)$ is smooth and flat at $x=0$, 
then the centralizer 
 $Z_\varphi$ coincides with the 1-parameter family 
$\{\exp (tX_\varphi)=\varphi^t \,;\, t\in \R\}$   
generated by $X_\varphi$. 

In general case, using the order of real numbers,  
the centralizer $Z_\varphi$ turns out to be a totally 
ordered abelian group which contains $\{\varphi^\Z\}\cong\Z$.  
Therefore it is uniquely identified with a certain 
subgroup of the additive group $\R$ under the identification  
$\{\varphi^\Z\}\cong\Z$.  
Depending on $\varphi$, 
$Z_\varphi$ can be 
far beyond the normal expectation, 
\EG it can be 
$\Z$, $\Q$, or 
$\Z\oplus\Z\alpha$ where  
$\alpha$ is a Liouville number \cite{Ey}, 
or far more complicated. 
The topology on $Z_\varphi$ through this identification 
with natural topology of $\R$ coincides with 
the one induced from 
the $C^0$-topology on 
$\mathit{Diff}^\infty([0,\infty))$.  

We should also remark that for 
$\,\varphi\,$ which is infinitely tangent to the identity 
at $\,x=0\,$,   
so is any element of $\,Z_\varphi,$.

At present it is not known whether $Z_\varphi\cong\R$ 
implies the smoothness of $X_\varphi$ at $x=0$. 
This is a subtle point in Hilbert's problem No. 5 when 
it is stated in the context of a continuous 
homomorphism from 
a Lie group to a group of diffeomorphisms.  
The difficulty occurs when the orbit is not compact.  
\subsection{Structure of $\mathit{Aut} R$}
\label{structure}
Upon all the previous preparations 
we are able to describe the structure of 
$\mathit{Aut} R$ as follows. 
\begin{prop}{\rm \quad Let $R$ be a Reeb component 
of real dimension $3$ which is given 
by the Hopf construction.  
\\
1)\quad 
The group $\mathit{Aut} R$ 
of automorphisms of the Reeb component $R$ 
admits a following 
sequence of extensions by abelian groups 
$\,\mathit{Aut}_0 H\,$, $\,Z_\varphi\,$, 
and $\,\Ker\,$,  
$$
0 \to \mathit{Aut}(R, H) \to \mathit{Aut} R 
\to \mathit{Aut}_0 H \to  0 
$$
$$0 \to \Ker \to \mathit{Aut}(R, H) 
\to Z_\varphi \to  0 
\\
$$
where $\mathit{Aut}_0 H \cong \C^*/\lambda^\Z$ is 
the multiplication by the constant linear part $a$ 
mod $\,\lambda^\Z\,$ as 
described in Lemma \ref{principal lemma}, 
$Z_\varphi$ is the centralizer of $\varphi$ 
which is explained in the previous subsection, 
and 
$\Ker$ 
is identified with the space of solutions to 
the equation (b) in Lemma \ref{principal lemma}.  
As explained in Part II, 
$\,\Ker\,$  is isomorphic to an infinite dimensional 
vector space $\solZ$ in Case (2) and (3), 
while in Case (1) it is a complex vector space 
of dimension 1 or 0 according to 
the resonance condition $\,\lambda=(\varphi'(0))^n$ 
for some $\,n\in\N\,$ or not.   
}
\end{prop}
The following is an important consequence to 
Theomre \ref{Main} 
in Part II on $\,\Ker\cong\solZ \,$ in Case (2) and (3).   
\begin{cor}\label{extension}{\rm
\quad  
If we paste $H\times (-\varepsilon, 0]$ to $R$ 
along the boundary $H$, 
in Case (3) any element of 
$\mathit{Aut}(R, H)$ 
extends to the other side, being the identity on 
 $H\times (-\varepsilon, 0]$, 
as a diffeomorphisms of class $C^\infty$.  
In Case (2) the same applies to $\,\Ker\,$.  
}
\end{cor}
\noindent
{\it Proof} of Proposition. \quad 
The first step of the extensions is obtained 
by looking at the action on the boundary, 
and once we assume that 
the action on the boundary is trivial, 
the second extension is obtained by looking at the 
action on the vertical line $\{0\}\times[0,\infty)$.  
We can interpret it as an action on the leaf space. 
\qed

In the above, 
the first extension does not yield a non-abelian group. 
Using the identification 
$\mathit{Aut} R\cong \mathit{Aut}(\tilde R; T)/T^\Z$ 
in Corollary \ref{identification}, 
we obtain a better description 
not only from the above point of view 
but also from that of the question 
whether the restriction map 
$r_H:\mathit{Aut} R \to \mathit{Aut}_0 H$ 
admits a homomorphic section. 
Note that $Z_\varphi$ admits a section to 
$\mathit{Aut}(R, H) \subset \mathit{Aut} R$. 

An element $f \in \mathit{Aut}(\tilde R; T)/T^\Z$ 
admits a presentation $f(z,x)=(az+b(x), \eta(x))$ 
up to $T^\Z$ where $T(z,x)=(\lambda z, \varphi(x))$.   
Therefore ignoring $b(x)$ from this presentation 
and assigning $f \mapsto (a,\eta)$ (mod $(\lambda, \varphi)^\Z$), 
we obtain a surjective homomorphism 
$\mathit{Aut}(\tilde R; T)/T^\Z \twoheadrightarrow 
(\C^*\times Z_\varphi)/(\lambda, \varphi)^\Z$ to an abelian group.  
Also, by setting $b(x)=0$, we see this abelian group 
can be realized as a subgroup of $\mathit{Aut}(\tilde R; T)/T^\Z$. 
This enables us to describe the structure 
of $\mathit{Aut} R$ as follows.  
\begin{thm}\label{MainThm}{\rm
\quad
The automorphism group 
$\mathit{Aut} R\cong \mathit{Aut}(\tilde R; T)/T^\Z$ 
is isomorphic to the semi-direct product 
$$
\Ker\, 
\rtimes \,\left\{
(\C^*\times Z_\varphi)/(\lambda, \varphi)^\Z
\right\}
$$
where $a\in\C^*$ acts on $b(x)\in \Ker$ 
by multiplication 
$b(x)\mapsto a^{-1}b(x)$, \IE the conjugation in 
the affine transformations of each leaf,   
and 
$\eta\in  Z_\varphi$ acts 
by $b(x)\mapsto b(\eta(x))$.   
}
\end{thm}
{\it proof.}\quad 
Let us only verify the action of $a$.  
The conjugation by the multiplication by $a$ 
is $[z\mapsto z+b(x)] \mapsto 
[z \mapsto a^{-1}(az + b(x))=z+a^{-1}b(x)]$. 
\hspace*{\fill} $\square$
\\

To close this section, consider the liftability of 
$\mathit{Aut}_0 H$ 
to $\mathit{Aut} R$.  
This is  nothing but the liftability of the surjective 
homomorphism 
$$
(\C^*\times Z_\varphi)/(\lambda, \varphi)^\Z 
\twoheadrightarrow 
\C^*/\lambda^\Z \, .
$$
Here we assume the continuity of splitting, 
otherwise the question should include thinking about 
non-continuous homomorphism 
$\R \to \R$ with $1\mapsto 1$.  
If the centralizer $Z_\varphi$ is the total of $\R$, 
it implies $Z_\varphi$ 
is a 
$C^0$-family 
of 1-parameter subgroup 
$\{\eta_t\, ;\, t\in\R\}$ 
in $\mathit{Diff}^\infty([0,\infty))$ 
with $\varphi = \eta_1$. 
Then we obtain easily a lift 
defined as 
$$
a (\mathrm{mod}\, \lambda^\Z ) 
\mapsto 
(a, \eta_{t(a)})\,
(\mathrm{mod}\, (\lambda, \varphi)^\Z )\, ,
\quad 
t(a)=\frac{\log \vert a\vert}
{\log\vert\lambda\vert} 
\, .
$$ 
The converse is almost the same. 
If we have a continuous lift 
to $\mathit{Aut}(\tilde R; T)/T^\Z$, 
choose a value of $\log \lambda$ and 
look at the lift of a circle subgroup 
$e^{t\log \lambda}$ ($0\leq t\leq 1$) to 
a continuous path in 
$\mathit{Aut}(\tilde R; T)$ 
starting from the identity. 
Then its projection to $Z_\varphi$ 
gives rise to a 1-parameter family 
in $Z_\varphi$ starting from the identity 
which ends at 
$\varphi$.  
If this curve is smooth, it implies that 
the Szekeres vector field $X_\varphi$ of 
$\varphi$ is smooth. 
Thus we obtain the following. 
\begin{thm}\label{splitting}{\rm
\quad  The restriction map 
$r_H:\mathit{Aut} R \to \mathit{Aut}_0 H$ admits 
a continuous [{\it resp.} smooth] homomorphic section 
if and only if the centralizer $Z_\varphi$ 
is isomorphic to $\R$ as an ordered abelian group 
[{\it resp.} the Szekeres vector field $X_\varphi$ is smooth]. 
In such cases, 
$\mathit{Aut} R$ admits not only 
a structure of twice semi-direct products 
$$
\mathit{Aut} R \cong 
(\Ker \rtimes Z_\varphi) 
\rtimes \mathit{Aut}_0 H
\cong 
(\Ker \rtimes \R) 
\rtimes \R^2/\Z^2 \, ,
$$
but also a structure of simple semi-direct product 
$$
\mathit{Aut} R \cong 
\Ker \rtimes 
(\mathit{Aut}_0 H \times Z_\varphi) 
\cong 
\Ker \rtimes 
(\R^2/\Z^2  \times \R). 
$$
of two abelian groups, 
where the action of the right group on the left 
is continuous with respect to the smooth topology 
on $\Ker\cong\solZ \subset C^\infty([0, \infty);\C)$ 
in Case (2) and (3) 
and $\Ker\cong\C$ or $\{0\}$ in Case (1) 
[{\it resp}. smooth in a usual sense].  
}
\end{thm}

\begin{rem}{\rm \quad 
We saw that in Case (2) and (3) 
the autmorphism groups are of infinite  
dimension, while in Case (1) it is of finite 
dimension and shows a similarity to 
compact complex manifold.

In fact, in Case (3), not only any Reeb component 
is realized in a Hopf surface 
in Construction \ref{HopfSurface},  
but also any of its automorphims 
extends to the ambient Hopf surface. 
Undr the setting $n=1$, the Hopf surface $W$ 
is obtained as $W=(\C^2\setminus \{O\})/T''$ 
where $T''(z, w)=(\lambda\cdot z,\, \mu\cdot w)$.  
Then a Levi-flat hypersurface $M= (\C\times\R \setminus\{O\})/T''$ 
contains a Reeb component $R$ of Case (3) 
for $\lambda$ and $\mu\,$. It is easy to see that 
the automorpshim 
$(z, w)\mapsto (a\cdot z + b w^p\, , \, c\cdot w)\,$
of the Hopf surface $W$ 
restricts to $R$ 
where $a$, $c \in  \C^*$ and $b\in\C$ are 
arbitrary constants in the resonant case $\lambda=\mu^p$ and 
$b=0$ in the non-resonant case. 
Therefore any 
of $\mathit{Aut} R$ extends to $W$.

On the other hand, 
a recent work of Koike and  Ogawa \cite{KO} 
seems suggesting that Reeb componets of Case (3) 
never appers in a Levi-flat hypersurfaces 
in a complex surface.  
Our result also mildly suggests that 
the same might apply to Case (2).   

Even in Case (1) it should be still confirmed 
whether if the automorphism extends to the ambient surface 
in the case where the Reeb component appears 
as a part of a Levi-flat hypersurface which bounds 
a Stein surface.   
}
\end{rem}

\section{Reeb foliations}
The automorphism group 
of a leafwise complex foliation 
on a closed 3-manifold which consists of 
two Reeb components is now easy to compute.  
In this section we assume 
that the relevant holonomy is infinitely 
tangent to the identity at the origin, 
because it must be so 
except for two special cases where 
two Reeb components are pasted in the same direction 
of the holonomies or exactly in the inverse direction, 
in both of which cases the pasting yields 
$\,S^2 \times S^1\,$.   

Let $R_{\varphi, \lambda}$ be the Reeb component 
which we dealt with in the previous section. 
For another pair of a diffeomorphism 
$\psi \in \mathit{Diff}^\infty([0,\infty))$ 
which is also expanding and infinitely tangent to the identity 
at the origin 
and a constant 
$\mu\in\C$ with $\vert\mu\vert>1$,  take the Reeb component 
${R}_{\psi, \mu}$ and let $\overline{R}_{\psi, \mu}$ 
denote the mirror of ${R}_{\psi, \mu}$, namely 
the one which we obtain by reversing the 
the transverse orientation. 
It is done by replacing  $x$ and $\varphi(x)$ 
with $-x$ and $-\varphi(-x)$ in the Hopf construction. 

For example if $\lambda=\mu$ we can paste  
$R_{\varphi, \lambda}$ and $\overline{R}_{\psi, \mu}$ 
along the common boundary $H=\C^*/\lambda^\Z$ 
by the identity of $H$ 
to obtain a leafwise complex foliation on $S^2\times S^1$. 
In general according to the pasting element 
$\in \mathit{SL}(2;\Z)$ we can choose appropriately 
$\lambda$ and $\mu$ and paste them. 
The foliation on $S^3$ obtained in such a way is called 
the Reeb foliation. 

Corollary \ref{extension} yields the following results.   
\begin{thm}{\rm 
\quad 
Let $(M, \F, J)$ be a leafwise complex foliation 
which is obtained by pasting  
$R_{\varphi, \lambda}$ and $\overline{R}_{\psi, \mu}$.  
Then its group of automorphism is 
naturally isomorphic to the fibre product 
of $\mathit{Aut} R_{\varphi, \lambda}$ 
and 
$\mathit{Aut} \overline{R}_{\psi, \mu}$ 
with respect to $\mathit{Aut}_0 H$.  

If the centralizer $Z_\psi$ is isomorphic to 
$\R$ as an ordered abelian group, then 
 $\mathit{Aut} R_{\varphi, \lambda}$ is continuously 
realized as a subgroup in the resulting group of automorphisms. 
}
\end{thm}
\begin{thm}{\rm 
\quad
If $(M^3, \F_1, J_1)$ is obtained from 
$(M^3, \F_0, J_0)$ by turbulization along 
a closed transversal and the resulting Reeb component 
is isomorphic to  $\mathit{Aut} R_{\varphi, \lambda}$, 
the group  $\mathit{Aut}(M^3, \F_1, J_1)$ naturally contains 
a subgroup which is isomorphic to 
 $\mathit{Aut} (R_{\varphi, \lambda}, H)$.  
}
\end{thm}

\begin{rem}{\rm\quad 
In both of above theorems, the automorphism group 
contains an infinite dimensional vector space $\mathcal Z$ 
or one more copy.  
Thus even in the case of closed manifolds, 
the automorphism group of leafwise complex foliation 
can be fairy large. 
This  presents a clear contrast between 
a leafwise complex foliation on a compact manifold 
and a complex structure on a compact manifold.   
Due to Meersseman and Verjovsky \cite{MV} 
in the study of moduli spaces, 
they present similar features 
as far as we deal with 
tame leafwise complex foliation.   
}
\end{rem}
\vspace{30pt}

\noindent
{\bf \Large  Part II : Schr\"oder's equation on the half line}
\\

\noindent
We study the functional equations on the half line $[0, \infty)$ 
which appeared in Section 3. 
The simpler one takes the form
$$
\beta\circ \varphi (x) = \lambda \beta(x)
$$
for a fixed diffeomorphism 
$\varphi$ and a constant $\lambda$. 

Ernst Schr\"oder 
started to study 
a similar (in fact, formally the same) 
functional equation on the unit disk $D$ 
in the complex plane $\C$ 
under the complex analytic setting  
in \cite{Sch} in the late 19th century.  
Not only because 
it is just the natural eigenvalue problem 
for the pull-back operator  
to look for $\beta$ and $\lambda$ 
for a given $\varphi$, 
also Schr\"oder initiated complex dynamical studies 
and was interested in the iteration of compositions 
maybe in the context of Newton's method.  
According to the development of the complex dynamics 
the problems that were treated in these epoch 
has become fairy well-understood.  
Recently the studies in this direction seem to be  
aiming at higher dimensional cases.  
For the history of an early stage of the complex dynamics 
and Schr\"oder's functional equation, 
we have two nice references [Al1, 2].  

Our aim is to solve the equations 
\begin{eqnarray*}
{\mathrm{Equation \,\, (I)\,}} &:& \beta\circ \varphi (x) 
= \lambda \beta(x)
\\
{\mathrm{Equation \, (II)}} &:& 
\beta_1\circ \varphi (x) = \lambda \beta_1(x), 
\quad  
\beta_2\circ \varphi (x) 
= \lambda \beta_2(x) + \beta_1(x)
\end{eqnarray*}
on the half line $[0, \infty)$ 
for an expanding diffeomorphism 
$\varphi \in \Diff^\infty([0, \infty))$ 
and a complex constant $\lambda$ with $|\lambda | >1$.  
Equation (II) is generalized to still higher dimensional case 
(II') which is expressed by using vector notations as 
$$
\qquad \!\!\!
{\mathrm{Equation\,(II')}}\,\,\, :\, \, \,
\boldsymbol\beta\circ \varphi(x) 
=A \boldsymbol\beta(x) \qquad \qquad \qquad \qquad \qquad \qquad
\qquad \quad 
$$
where $\varphi$ is as above, 
$A=(a_{ij})$ is an $M\times M$ matrix 
any of its eigenvalues 
has the absolute value greater than $1$,  
and 
$\boldsymbol\beta(x)
= {}^t(\beta_1(x), \cdots ,  \beta_M(x)) 
\in C^\infty([0, \infty); \C^N)$ 
is the unknown function.  
Of course, the problem is easily reduced to the case 
where $A$ is a single non-trivial Jordan block.  
So we can assume that $A$ has a unique eigenvalue 
$\lambda$ as above and eventually  
$a_{ij}=\lambda$ for $i=j$, $a_{ij}=1$ for $i=j+1$, 
and $a_{ij}=0$ otherwise. 

It is easily seen that in the case of $\varphi'(0)=1$, 
which is of our main concern, 
there is no analytic solution but 
$\beta(x)$ or $\boldsymbol\beta(x)\equiv 0$. 
On the other hand, if  $\varphi'(0)>1$, 
for Equation (I) the space of solution is trivial or 
of dimension 1, depending on the resonance of 
$\lambda$ and $\varphi'(0)$.   
These are in fact exactly the same even when 
working on $D\in \C$.  
Therefore Schr\"oder's equation 
exhibits very distinctive feature 
when it is considered on the half line 
with $\varphi'(0)=1$.

The space of solutions to Equation (I) 
turns out to be an infinite dimensional vector space 
which is in a sense 
isomorphic to $C^\infty(S^1; \C)$ 
whenever $\varphi'(0)=1$.  
In the subsequent sections, 
first we describe the space of solutions much clearer, 
and then we give two different proofs.   

It is to be remarked that as  $\varphi$ is expanding, 
our problem is essentially that on the germs around $x=0$. 
For a given germ of $\varphi$, 
we can extend $\varphi$ to the whole of $[0, \infty)$ 
as a realization of the germ as far as $\,\,\varphi(x)>x\,\,$ 
is satisfied for $\,x>0\,$.  
Then the same applies to $\beta(x)$ 
because once it is given as a germ around  $x=0$, 
it is automatically and uniquely extended to the whole half line 
by the equation itself.

The results on Equation (I) 
are used in Part I of this paper. 
Those for (II) serve in \cite{Ho} to determine 
the automophism groups of Reeb components 
with complex leaves of complex dimension 2.  
When we extend our results on the automorphism groups 
to higher dimensional cases, 
the results on (II') are necessary.

\section{\large The space of solutions}
In this section we give precise statements of our results 
on the Schr\"oder type functional equations (I), (II), and (II') 
and describe the spaces of their smooth solutions. 
\\

For an expanding diffeomorphism 
$\varphi\in\mathit{Diff}^\infty([0, \infty))$ 
and a complex number $\lambda$ with $\vert\lambda\vert > 1$, 
we consider Equation (I), (II), and (II').  
First consider these equations 
(not on the whole $[0, \infty)$ but) on $(0,\infty)$. 
Then, Equation (I) has a lot of solutions 
and if we fix any solution 
$\beta^*(x)\in C^\infty((0,\infty); \C)$ 
which never vanishes, \IE $\beta^*(x)\ne 0$ 
for $x>0$, 
then each solution corresponds to a 
smooth function on $S^1=(0,\infty)/\varphi^\Z$
by assigning $\beta \mapsto \beta/\beta^*$. 
This gives a bijective correspondence between 
the space ${\mathcal Z}=\solZ$ of solutions  to (I) 
on $(0, \infty)$ 
and $C^\infty(S^1; \C)$ as vector space. 
This correspondence will be more precise 
in  Section \ref{Fourier}  
in fixing the coordinate on the circle.  

In Section \ref{Fourier}, 
we will make this correspondence 
more precise, 
fixing 
the coordinate on $\,S^1\,$ and the choice of 
$\,\beta^*\,$.

Also take the space ${\mathcal S}=\solS$ of solutions to Equation (II) 
on  $(0, \infty)$. If we assign $\beta_1$ to  a solution 
$(\beta_1, \beta_2) \in \mathcal S$, we obtain the projection 
$P_1: \mathcal S \to \mathcal Z$. 
Here the kernel 
of $P_2$ is nothing but $\mathcal Z$. 
We also see that the projection $P_1$ is surjective because 
for any $\beta_1\in \mathcal Z$
$$
\beta_2(x)= \frac{1}{\lambda\log\lambda} \beta_1(x)\log\beta^*(x)
$$
gives a solution $(\beta_1, \beta_2)\in \mathcal S$, 
where for $\log\lambda$ and $\log\beta^*(x)$ 
any (smooth) branch can be taken. 
Therefore, as a vector space,  
$\mathcal S$ has a structure such that 
$$
0 \to \mathcal Z \to \mathcal S \to \mathcal Z \to 0
$$
is a short exact sequence. 

For Equation (II') we simply repeat this extension.  
Let 
${\mathcal S}_M$ 
denote the set of solutions on $(0, \infty)$. 
For 
$1\leq m<M$, 
 ${\mathcal S}_m$ 
is naturally identified with 
a quotient 
$\{{}^t(\beta_1,\, \cdots , \, \beta_m) \,\vert\,    
{}^t(\beta_1, \, \cdots ,\, \beta_M) \in  {\mathcal S}_M
\}$ 
of  ${\mathcal S}_n$. 
Each projection 
$
{\mathcal S}_m \to {\mathcal S}_{m-1}
$ 
is surjective because 
the multiplication 
$
\times \frac{1}{\lambda\log\lambda}\beta^*
$
is a linear right inverse 
and its kernel coincides with $\mathcal Z$. 
${\mathcal S}_1$ is nothing but $\mathcal Z$ and 
${\mathcal S}_2$  the above ${\mathcal S}$ as well.  

Let $\,\varphi^*\,$ denote the pull-back by $\,\varphi\,$. 
Then Equation (I) is expressed as 
$\,\,(\varphi^* - \lambda)\beta=0\,\,$.   
Now, Equation (II) is nothing but 
$\,\,(\varphi^* - \lambda)^2\beta_2=0\,$, 
where we put  $\,\,\beta_1\,\,$ 
by setting 
$\,\,(\varphi^* - \lambda)\beta_2=\beta_1\,$. 
Inductively we see Equation (II') 
is nothing but 
$\,\,(\varphi^* - \lambda)^M\beta_M=0\,\,$ while 
$\,\,(\varphi^* - \lambda)\beta_m=\beta_{m-1}\,\,$ 
for $\,m=2,3,4, \cdots, M\,$ might also be regarded as 
auxiliary equations.   \vspace{8pt}

As in Section 3 of Part I, 
we devide the situation into the following three cases 
according to the nature of the jet of $\varphi$ at $x=0$.  
For the second and the third cases, 
the statements of our results are the same.   
Here, $f^{(i)}$ denotes the $i$-th derivative of $f(x)$ 
and $j^if(0)$ denotes the $i$-th jet at $x=0$. 
\vspace{5pt}
\\
\noindent
\quad Case (1) : $\varphi'(0)=\mu >1$ , 
\vspace{2pt}
\\
\noindent
\quad Case (2) : for some $n \geq 2$ ,   
$j^{n-1}\varphi(0)=j^{n-1}\mathit{id}(0)\,$  
and 
$\,\,\varphi^{(n)}(0)>0$ ,  
\vspace{2pt}
\\
\noindent
\quad Case (3) : 
$j^\infty\varphi(0)=j^\infty\mathit{id}(0)\,$ .  
\vspace{8pt}

Let us state the results. We start with the easiest case.

\begin{thm}\label{linear expansion}{\rm
\quad Consider Case (1). 
\\
1) (Resonant case)\quad If $\,\,\lambda=\mu^n\,\,$ 
is satisfied for some $\,n\in\N$ ,    
then  the space of solutions  $\K$ 
to Equation (I) 
is a complex vector space of dimension 1. 

For Equation (II'), $\,\boldsymbol\beta\,\,$ satisfies  
(II') if and only if 
$\,\,\beta_1= \cdots =\beta_{n-1}\equiv 0\,\,$ and 
$\,\,\beta_n\in\K\,\,$ hold. 
Therefore in total the space of solution is also 
1-dimensional. 
\vspace{3pt} 
\\
2)(Non-resonant case)\quad If no positive integer 
$n\in\N$ satisfies  $\lambda=\mu^n$, 
then there exists no solution to (I) but 
$\beta(x)\equiv 0$, and we have $\K=\{0\}$.   

For Equation (II') the same applies. 
Namely the only smooth solutions on $\,[0,\infty)\,$ 
is $\,\,\boldsymbol\beta(x)\equiv 0$ .  
}\end{thm}
Remark here that if $\varphi(x)=\mu x$, 
the solution to (I) in resonant case 
is nothing but $\beta(x)=\mathrm{const}\cdot x^n$ .  
Also remark that accordingly $\lambda$ is a 
positive real number. 
This result is so easy that the proof is given here. 
\vspace{8pt}
\\
\noindent
{\it Proof.} \quad 
From Sternberg's linearization theorem \cite{St}, 
there exists a diffeomorphism 
$\,h\in\Diff^\infty([0, \infty))\, $ 
which conjugates $\varphi$ into the linear one  
$\, h^{-1}\circ\varphi\circ h(x)=\mu x\,$.  
Therefore in solving the equations,  
from the first we can assume 
$\,\varphi(x)=\mu x\,$ . 
Therefore the equation (I) takes the following form. 
$$
\beta(\mu x)=\lambda\beta(x)\qquad 
\mathrm{for}\,\,\, x\in [0,\infty)\, .  
$$
In both cases, 
by differentiating Equation (I) for arbitrary many times 
at $x=0$, we see that the Taylor expansion at $x=0$ can be 
non-trivial only at the degree $n=\log\lambda/\log\mu$.  
Therefore in the resonant case, 
$\beta(x)$ must be in the form 
$\beta(x)=c\cdot x^n +f(x)$ where $f(x)$ is a flat function, 
and in the non-resonant case, the same form with $c=0\,$. 
Then, in the resonant case, 
as $c\cdot x^n$ is a solution to (I), so is $f(x)\,$.  
Therefore in both cases, it is enough to check that 
any non-trivial flat function can no be a solution. 

If we had a non-trivial flat 
solution $f(x)\,$, 
it would contradict as follows. 
Take $x_0 \in (0, \infty)$ with $f(x_0)\ne 0$ and 
look at $f(\mu^{-k}x_0)=\lambda^{-k}f(x_0)$ for $k\in\N$.  
On the other hand, as $f$ is flat we have 
$\lim_{x\to 0} f(x)/x^l =0$ for any $l\in \N$. 
So large enough $l$ ($\geq |\log\lambda/\log\mu|$) 
gives rise to a contradiction.  

This is well-known also as a fact (even for 
higher dimensional case) 
that a weighted-homogeneous function 
is smooth at the origin only when 
it is a polynomial. 

For Equation (II), from the above result, 
we assume $\beta_1(x)=cx^n$ in the resonant case.  
Then a similar computation for $\beta_2$ 
implies $c=0$. Therefore we have $\beta_1=0$ and thus 
$\beta_2=c'x^n$ for some $c'\in \C$.  
In the non-resonant case, 
the argument for (I) suffices. 
For (II'), the argument for (II) works as an inductive step. 
\qed

The results for Case (2) and (3) 
can be stated together.

\begin{thm}\label{Main}\quad{\rm
For both of Case (2) and Case (3), 
and for any $\,\,\lambda\in \C\,\,$ with $\,\,|\lambda|>1\,$,  
the followings hold.  \vspace{3pt} 
\\
1) Any solution $\,\,\beta\in\mathcal Z_{\lambda, \varphi}\,$ 
to Equation (I) 
on $\,\,(0, \infty)\,\,$ 
extends to $\,\,[0,\infty)\,\,$ 
so as to be a smooth function which is flat at $\,x=0\,$, 
\IE the $k$-th jet satisfies 
$\,\,j^k\beta(0)=0\,\,$ for any $\,k=0,1,2, \,\cdots$. 
In other words, the space $\,\,\K\,\,$ of all solutions to (I) 
considered on $\,\,[0,\infty)\,\,$ 
coincides with $\,\,\mathcal Z_{\lambda, \varphi}\,\,$. 
\vspace{3pt} 
\\
2) The same applies to Equation (II) and (II').  
Namely by putting $\,{\boldsymbol\beta}(0)=0\,$, 
any $\,\,{\boldsymbol\beta}\in {\mathcal S}_M\,\,$ 
is smooth and flat at $\,x=0\,$. 
}
\end{thm}

In the next section we give a proof for Case (3), 
which is more direct and simpler than 
one given in the final section. 
In this proof, 
we directly estimate the derivatives of any order 
of $\beta\circ\varphi$ for $\beta \in \mathcal Z$.   
The higher derivatives of a composite function 
is complicated and described in the formula 
of Fa\`a di Bruno.  We do not need its full length.  

Unfortunately this method does not work for Case (2). 

A proof by a different approach for Case (2) is 
given in the following section. 
It is based on two theories.  
One is Takens' normal forms \cite{Ta} 
for germs around the origin 
of $\,\Diff^\infty([0,\infty))\,$ 
and of vector fields on $[0, \infty)$  
with non-trivial finite order jets.   
The other one is a classical theory of 
Fourier expansion of $C^\infty(S^1; \C)$.  
Unfortunately again, this method 
seems difficult to apply to Case (3).  

Then in the final section we give a proof 
relying on the center manifold theory,  
which covers both of Case (2) and (3).   
This proof is suggested by Masayuki Asaoka.  
It might be worth remarking that 
when a proposition is proved in the framework of 
hyperbolic dynamical systems, 
quite often it is also proved in Fourier analysis, 
and vice versa.  
Here we might observe a similar phenomenon.

\section{\large Direct proof for Case (3)}
We prove Theorem \ref{Main} for Case (3), 
namely, 
in the case where $\varphi$ is flat to the identity, 
by a direct estimate of the derivatives of $\beta(x)$ 
of an arbitrary order when $x \to 0$. 

In order to clarify the strategy 
it might be suggested to the readers to check 
$\ds\lim_{x\to +0}\beta(x)=0$ 
and 
$\ds\lim_{x\to +0}\beta'(x)=0\,$, 
\IE for $k=0, 1$,  
which are easy and reviewed in Proposition \ref{1jet}, 
and then the the second jet $k=2$). 
Looking at up to the case $k=3$ might make 
the roll of the following lemma clearer.

\begin{lem}\label{Faa di Bruno}\quad{\rm 
The $n$-th derivative $\{\beta(\varphi(x))\}^{(n)}$ 
is written in the following form for $n\in\N$. 
$$
\{\beta(\varphi(x))\}^{(n)}
=(\varphi'(x))^n \cdot\beta^{(n)}(\varphi(x))
+
\sum_{k=1}^{n-1}\Phi_{n,k}\cdot\beta^{(k)}(\varphi(x))\, .
$$
Here, $\Phi_{n,k}$ is an integral polynomial in 
$\varphi'(x)$, $\varphi''(x)$,$\cdots$,  $\varphi^{(n)}(x)$, 
without constant term and no term is of monomial 
only in $\varphi'(x)$. 
}
\end{lem}

This lemma is easily seen by the induction,  
but in fact it is a corollary to 
the well-known formula of Fa\`a di Bruno 
(\EG see \cite{Ri}, \cite{Ro}, or textbooks on calculus). 
It is independent of our assumption on $\varphi$ and 
is valid for any composite functions.    
On the other hand the flatness of $\varphi$ at the origin 
implies  
$(\varphi'(x))^n \to 1$ and $\Phi_{n,k}\to 0$ 
when  $x \to 0+0$. 
\\

Now let us prove 1) of Theorem \ref{Main}. Let $\beta$  be 
a solution to (I) on $(0,\infty)$.  
From the equation it is easy to see that $\beta(x)\to 0$ 
when $x \to 0+0$.  

Now fix any integer $N$. 
 $\beta'(x)\to 0$ is also easy to see, but for higher 
derivatives, in a natural estimate the lower derivatives 
are involved. Thus the basic strategy is not to estimate 
the higher derivatives by induction on the order, 
but to estimate them all together up to the 
fixed order $N$. 

From Equation (I) and the above lemma 
we have the following computation. 
\begin{eqnarray*}
\sum_{n=1}^{N}\vert\beta^{(n)}(x)\vert
 &=& 
\frac{1}{\vert\lambda\vert}
\sum_{n=1}^{N}
\left\vert 
\{\beta(\varphi(x))\}^{(n)}
\right\vert 
\\
 &\leq& 
\frac{1}{\vert\lambda\vert}
\sum_{n=1}^{N}
\left\{
(\varphi'(x))^n \cdot\vert\beta^{(n)}(\varphi(x))\vert
+
\sum_{k=1}^{n-1}\vert\Phi_{n,k}\vert
\cdot
\vert\beta^{(k)}(\varphi(x))\vert
\right\}
\\
 &\leq& 
\frac{1}{\vert\lambda\vert}
\sum_{k=1}^{N}
\left((\varphi'(x))^k +
\sum_{n=k+1}^{N}\vert\Phi_{n,k}
\vert\right)
\cdot
\vert\beta^{(k)}(\varphi(x))\vert
\end{eqnarray*}
As is remarked above, we know $(\varphi'(x))^k \to 1$ 
and 
$\sum_{n=k+1}^{N}\vert\Phi_{n,k}
\vert \to 0$ when $x\to 0$. Therefore there exists 
$b_N>0$ such that for $x\in (0,b_N]$ 
we have   
$$(\varphi'(x))^k +
\sum_{n=k+1}^{N}\vert\Phi_{n,k}
\vert \leq \sqrt{\vert\lambda\vert} 
\quad 
\mathrm{for}\,\,\,\, k=1, 2, \cdots , N .
$$
This implies for any $x \in (0, b_N]$ 
$$
\sum_{n=1}^{N} \vert \beta^{(n)}(x) \vert 
\leq 
\frac{1}{\sqrt{\vert\lambda\vert}}
\sum_{n=1}^{N} \vert \beta^{(n)}(\varphi(x)) \vert .
$$

Put
$
M=\max
\{
\sum_{n=1}^{N} \vert \beta^{(n)}(x) \vert 
\,; \, 
x\in[b_N, \varphi(b_N)]
\}
$ 
and define $m(x)\in\N$ for $x\in (0, b_N)$ 
so that 
$\varphi^{m(x)}\in[b_N, \varphi(b_N))$. 
Then, the above inequality implies 
$$ 
\sum_{n=1}^{N} \vert \beta^{(n)}(x) \vert 
\leq
M\cdot \sqrt{\vert\lambda\vert\,}^{-m(x)}
$$
for  $x\in (0, b_N)$.   
Because `$x\to 0+0$' is equivalent to `$m(x)\to\infty$', 
we obtained the convergence 
$$
\beta^{(n)}(x) \to 0 \quad (x\to 0+0)
\qquad \mathrm{for} \quad n=1,\cdots, N .
$$
This completes the proof of 1). 
\\

Let us outline the proof of 2) for $M=2$. 
We extend the basic strategy of the proof of 1) 
in the following sense. When we estimate the derivatives 
of $\beta_2$, naturally those of $\beta_1$ are involved. 
Therefore we  will estimate the derivatives of $\beta_2$ 
and  $\beta_1$ all together up to a fixed order $N$, 
even though the flatness of $\beta_1$ is already proved in 1).  

First we fix 
$\varepsilon>0$ so small that 
$\displaystyle
\varepsilon \leq 
\vert\lambda\vert^{\frac{5}{4}}-
\vert\lambda\vert
$ 
is satisfied.  
Now take any solution $(\beta_1, \beta_2) \in \mathcal S$ 
and put $\tilde\beta_1=\varepsilon^{-1}\beta_1$. 
Then instead of Equation (II), $\tilde\beta_1$ and  $\beta_2$ satisfy 
\vspace{5pt}
\\
$
\quad \mathrm{Equation\,\, (\tilde{II})} :\quad 
\tilde\beta_1(\varphi(x))=\lambda \tilde\beta_1(x)
, \;\;\;  
\beta_2(\varphi(x))
=\lambda \beta_2(x)+ \varepsilon\tilde\beta_1(x)
\, . 
$
\vspace{5pt}
\\
Then, from ($\tilde{\mathrm{II}}$) we have  
\begin{eqnarray*}
\frac{e^{i\theta}\lambda - \varepsilon}
{\lambda}\tilde\beta_1(\varphi(x))
+
\beta_2(\varphi(x))
&=&
e^{i\theta}\lambda\tilde\beta_1(x)
+
\lambda\beta_2(x)
\end{eqnarray*}
and consider the $n$-th derivatives of both sides. 
For any $\theta\in\R$ and $n=1, \cdots , N$, we have 
$$
\vert e^{i\theta}\tilde\beta_1^{(n)}(x) + \beta_2^{(n)}(x)\vert
\leq 
\frac{1}{\vert\lambda\vert}
\left(
\frac{\vert\lambda\vert + \varepsilon}{\vert\lambda\vert} 
\vert \{\tilde\beta_1(\varphi(x))\}^{(n)}\vert 
+ 
\vert \{\beta_2(\varphi(x))\}^{(n)}\vert\right) 
$$
Because the right hand side is independent of $\theta$, 
using the inequality 
$$
\left\vert 
\frac{e^{i\theta}\lambda - \varepsilon}{\lambda}
\right\vert 
\leq
\frac{\vert\lambda\vert + \varepsilon}{\vert\lambda\vert} 
\leq 
\vert\lambda\vert^\frac{1}{4}
\quad
\mathrm{for\,\,\, any}\,\,\, \theta\in\R
$$ 
we obtain 
$$
\vert \tilde\beta_1^{(n)}(x)\vert + \vert\beta_2^{(n)}(x)\vert
\leq 
\frac{1}{\vert\lambda\vert^\frac{3}{4}}
\left(\vert \{\tilde\beta_1(\varphi(x))\}^{(n)}\vert 
+ 
\vert \{\beta_2(\varphi(x))\}^{(n)}\vert\right).   
$$

Applying Lemma \ref{Faa di Bruno} to $\tilde\beta_1(\varphi(x))$ and 
to $\beta_2(\varphi(x))$ 
for $n=1, \cdots , N$, 
from the same argument as in 1) we obtain 
$$
\sum_{n=1}^{N}
\left(\vert \tilde\beta_1^{(n)}(x)\vert + \vert \beta_2^{(n)}(x)\vert\right) 
\leq 
\frac{1}{\vert\lambda\vert^\frac{1}{4}}
\sum_{n=1}^{N}
\left(\vert \tilde\beta_1^{(n)}(\varphi(x))\vert 
+ \vert \beta_2^{(n)}(\varphi(x))\vert\right) 
$$
for $x\in (0, b_N]$, where $b_N$ is exactly the same as in the proof of 1). 
\hspace*{\fill}$\square$
\\

Now it is almost straight forward to further generalize 
this proof for Equation (II').

\section{\large Proof by 
Fourier series for Case (2)}\label{Fourier} 
A proof of Theorem \ref{Main} for Case (2) is given here. 
It relies on two big tools. 
The first one is 
Takens' normal form which plays a similar roll 
as Sternberg's linearization in the proof of 
Theorem \ref{linear expansion}.  
The second one is a classical Fourier expansion/series 
of smooth functions on the circle. 
Takens' normal form enables us to 
consider (countably many) simple linear homogeneous 
ordinary differential equations instead of 
considering Equation (I).  
Equation (II) and (II') correspond to inhomogeneous 
or vector valued case.  
As we will see below, 
we have an ODE for each choice of the value of 
$\,\log \lambda\,$ , 
Our functional equation (I) and 
countably many ODE's are related 
by Fourier expansion and series. 

\subsection{
Takens' normal form and Fourier basis
}
\begin{thm}\label{Takens}
{\rm(Takens, \cite{Ta})\quad 
Let $\,\varphi \in \Diff^\infty([0, \infty))\,$ 
be in Case (2).  \vspace{3pt}
\\
1)\quad 
There exists a diffeomorphism 
$h \in \mathit{Diff}^\infty([0,\infty))$ 
which conjugates $\varphi$ into a diffeomorphism 
$\psi\in \mathit{Diff}^\infty([0,\infty))$ 
of the following polynomial type 
on $[0, x_1]$ ($\exists x_1>0$) 
$$
\psi(x)=x + x^n + \alpha x^{2n-1} 
\quad \mathrm{and} \quad 
\psi = h^{-1}\circ\varphi\circ h \, .  
$$
The coefficient $\alpha\in\R$ 
is determined by the $(2n-1)$-jet of $\varphi$ at $x=0$. 
\vspace{3pt}
\\
\noindent
2) \quad 
There also exists a diffeomorphism 
$k \in \mathit{Diff}^\infty([0,\infty))$ which conjugates 
$\varphi$ 
in a neighborhood of the origin 
into 
the exponential map 
$$ 
\exp
X
= k^{-1}\circ\varphi\circ k \, .
$$
of the vector field 
$$
X=\rho(x) \frac{d}{dx}\, , 
\quad 
\rho(x)= x^n + a x^{2n-1} \, .  
$$ 
The the coefficients $\alpha$ in 1) and 
$a$ here are related by $a=\alpha - n/2\,$. 
}
\end{thm}
\begin{rem}
{\rm \quad 1) \quad 
The second statement follows from the first one, 
because a direct computation shows 
$$
\exp(x^n + a x^{2n-1})
\frac{d}{dx} = x + x^n + \left(a+\frac{n}{2}\right)x^{2n-1} 
+ [\mathrm{\,higher\,\, order\,\, terms}]\, . 
$$ 
2)\quad Takens also gave a normal form for vector fields.  
It takes almost the same form but we do not need it here. 
}
\end{rem}

Thanks to Takens' theorem, 
we can conjugate our equations by a smooth diffeomorphism 
and are allowed to assume that the holonomy 
$\varphi$ is of the form 
$$ 
\varphi=\exp X
\, , 
\quad 
 X
=\rho
(x)\frac{d}{dx}\, , 
\quad
\rho
(x)=x^n + a x^{2n-1}\,\, 
\mathrm{on} \,\,  
[0, x_0]
$$
for some $n\geq 2$, $a \in \R$, and $x_0>0\,$.  
We also assume that   
$\,\rho
(x)>0\,$ on $\,(0, \infty)\,$ and 
$\,\rho
(x)\equiv 1\,$ on $\,(x_1, \infty)\,$ 
for some  $x_1>x_0\,$.  
\vspace{8pt}

We consider the following ordinary differential equation on 
$(0, \infty)$    
\vspace{8pt}
\\
\qquad Equation (I-$\Lambda$) : \quad \qquad 
$\beta'(x)=\frac{\Lambda }{\rho(x)}\beta(x)\, . 
$
\vspace{8pt}
\\
\noindent
This is of course equivalent 
to the following ODE in the variable $t$. 
$$
\left.\frac{d}{dt}\right\vert_{t=0}
\beta(\exp(tX)(x))=\Lambda\cdot\beta(x) 
$$
Therefore any solution $\beta$ is presented as  
$\beta(\exp(tX)(x_0))=e^{\Lambda t}\cdot\beta(x_0)
=C\cdot e^{\Lambda t}$ 
for a constant $C\in\C$. 
It is also clear that $\beta$ satisfies 
Equation (I) on   
$(0, \infty)$. 
In the variable $x$, 
$\beta(x)$ is presented as 
$\beta(x)=C \cdot 
\exp\left( 
\Lambda \int_{x_0}^x \frac{1}{\rho(u)}
du\right)
$.  
In particular on $(0, x_0)$, we have 
$$\beta(x)=C \cdot 
\exp\left( 
(R + i\,{\mathrm{Im} \Lambda})
\int_{x_0}^x \frac{1}{u^n(1+au^{n-1})}
du
\right)
$$
%
%
where the real part 
$R=\log\vert\lambda\vert$ 
is positive.   

Now we choose $\Lambda_{\langle 0 \rangle}$ 
as a value of $\log \lambda$ 
and fix it.  
Then other general values of $\log\lambda$ 
are given as $\Lambda_{\langle l \rangle}= R +i\theta_l$,  
$\theta_l =\theta_0 + 2l\pi$ for $l\in \Z$. 
\vspace{8pt}

Here we make the correspondence between 
the space $\solZ$ of solutions of (I) considered on 
$(0, \infty)$ and $C^\infty(S^1;\C)$ more precise.  
As a coordinate on the circle $S^1$ we take 
$\theta = t$ (mod $2\pi$) where 
$x(t)=\beta(\exp(tX)(x_0))$ is assumed. 
Now for each $l\in \Z$, 
let $\beta_{\langle l \rangle}$ denote the solution to 
the ODE (I-$\Lambda_{\langle l \rangle}$) 
which satisfies $\beta_{\langle l \rangle}(x_0)=1$.   
Therefore we easily know that 
$\beta_{\langle l \rangle}(x(t))=
e^{2\pi l t \cdot i}\cdot 
\beta_{\langle 0 \rangle}(x(t))\,$.

Take $\beta_{\langle 0 \rangle}$ as $\beta^*$ 
in defining the correspondence. 
Then $\beta_{\langle 0 \rangle}$ corresponds 
to the constant function $1$ on $S^1$ and in general 
$\beta_{\langle l \rangle}$ 
corresponds to 
$\cbeta_{\langle l \rangle} \in  C^\infty(S^1;\C)\,$,  
namely,    
$$   
\cbeta_{\langle l \rangle}(\theta)=e^{2\pi l \theta \cdot i} 
\qquad \mathrm{for}\,\,\,  l \in \Z\, ,  
$$
so that 
$\cbeta_{\langle l \rangle}$'s ($l \in \Z$) 
form the standard Fourier basis 
for  $C^\infty(S^1; \C)$.  
The following is well-known and well fits into our situation. 
\begin{thm}{\rm (see \EG \cite{Ka})\quad 
The infinite sum with coefficients $c_l\in\C$ 
$$
\sum_{l=-\infty}^{\infty}c_l\cdot e^{i\theta}
$$
defines a smooth function on 
$\theta\in S^1=\R/2\pi\Z\,\,\,$ 
if and only if\,\, the sequence of coefficients 
$\{c_k\}_{k\in\Z}$ 
is rapidly decreasing, namely it satisfies 
$$
\sum_{l=-\infty}^{\infty}|l|^d|c_l|<\infty 
\quad \mathrm{for\,\, any}\,\,\, d\in\N 
\, .
$$
}\end{thm}
Therefore any $\beta\in\solZ$ is given as 
an infinite sum 
$$
\beta=\sum_{l=-\infty}^{\infty}c_l\cdot 
\beta_{\langle l \rangle}
$$
with a rapidly decreasing sequence of coefficients 
$\{c_k\}_{k\in\Z}$.

\subsection{Proof of Main Theorem}
What we have to prove in this section 
is stated as follows. 
\begin{thm}\label{restatement}
{\rm\quad 
For any rapidly decreasing 
$\{c_k\}_{k\in\Z}\,$,  
$\beta=\sum_{l=-\infty}^{\infty}c_l\cdot 
\beta_{\langle l \rangle}$ 
is extended to 
$[0,\infty)$ as $\beta(0)=0$ 
and is smooth and flat at $x=0$.  
}
\end{thm}
Let us verify this for each base.  
\begin{prop}{\rm\quad 
The solution  $\beta(x)$ to (I-$\Lambda$) is extended to 
$[0,\infty)$  as $\beta(0)=0$, 
and then  $\beta(x)$ is smooth and 
flat at $x=0$.  
}
\end{prop}

\noindent
{\it Proof.}\quad 
It is easy to compute the integration but 
we only need to remark 
that for some $\delta>0$ and any $x\in (0, \delta)$ we have 
$$ 
\left|\int_{x_0}^x \frac{1}{u^n(1+au^{n-1})}du\right| 
\geq 
\left|\frac{1}{2x}\right|\, . 
$$
The derivative of $\beta$ of order $k\in \N$  
is a multiplication of $\beta$ and some rational function in 
the variable $x$.  
Therefore for any $k \in \N$ we have 
$$
\vert\beta^{(k)}(x)\vert 
\leq 
\vert \mathrm{a\,\, rational\,\, function} \vert 
\times \exp(-\frac{1}{2x})
\, \to 0 \,\, (x\to 0)
$$ 
which suffices to show 
the smoothness and flatness of $\,\beta\,$ at $\,x=0\,$.  
\vspace{8pt}
\qed

In order to proceed further, 
we need to take 
a slightly closer look at those rational functions.   
Recall that $n$, $a$, and $\lambda$ are already fixed.  
\begin{lem}{\rm\quad 
For $k\in \N$ and $j=1, \cdots , k$, there exists a 
fixed polynomial $Q_{k,j}(x)$ which satisfies 
on $(0, x_0)$
$$
\beta_{\langle l \rangle}^{(k)}(x)=
\left\{
\frac{1}{P(x)^k}\sum_{j=1}^k Q_{k,j}(x)(R +i\theta_l)^j
\right\}
\beta_{\langle l \rangle}(x), 
\quad 
P(x)=x^n + a x^{2n-1}
$$
and $Q_{k,j}(x)$ is a linear combination of 
multiplications of $(k-j)$-many of 
$P(x)$, $P'(x)$, $\cdots$, $P^{(k-j)}(x)$, 
with total degree of differentiation $(k-j)$.   
}\end{lem}
For example,  
$Q_{k,k}(x)=1$,  $Q_{k,k-1}(x)=k(1-k)P'(x)/2$, and so on. 
The lemma is easily proved by the induction on $k$.  
\vspace{5pt}

Let us develop $(R +i\theta_l)^j$ into 
a polynomial of $l$ as follows. 
$$
(R +i\theta_l)^j
=
(R +i(\theta_0+2\pi l))^j
=\sum_{d=0}^jR_{j,d}l^d
$$
Here the constants $R_{j,d}$ 
($j\in\N$, $d=0$, $\cdots$, $j$) 
are determined by 
$R$ and $\theta_0$.  
\vspace{8pt}
\\

\noindent
{\it Proof} of Theorem \ref{restatement}. \quad 
For a rapidly decreasing sequence 
$$
\{c_k\}_{k\in\Z}
\quad 
\mathrm{with}
\quad
\sum_{l=-\infty}^{\infty}|l|^d|c_l| = M_d < \infty
\,\,\, 
\mathrm{for}\,\, \forall d \in \N\cup\{0\}   
$$
take 
$
\beta(x)=\sum_{l=-\infty}^{\infty}c_l\cdot 
\beta_{\langle l \rangle}(x)
$.    
Then we have the following estimate; 
\begin{eqnarray*}
|\beta^{(k)}(x)|
&=&
\left|
\sum_{l=-\infty}^{\infty}c_l\cdot
\beta_{\langle l \rangle}^{(k)}(x)
\right|
\\
&=&
\left|
\sum_{l=-\infty}^{\infty}c_l\cdot 
\left\{
\frac{1}{P(x)^k}\sum_{j=1}^k Q_{k,j}(x)(R +i\theta_l)^j
\right\}
\beta_{\langle l \rangle}(x)
\right|
\\
&\leq&
\left|
\frac{1}{P(x)^k}
\right|
\left|
\sum_{j=1}^k Q_{k,j}(x)
\sum_{l=-\infty}^{\infty}c_l
\left(
\sum_{d=0}^jR_{j,d}l^d
\right)
\beta_{\langle l \rangle}(x)
\right|
\\
&=&
\left|
\frac{1}{P(x)^k}
\right|
\left\{
\sum_{j=1}^k 
\sum_{d=0}^j
Q_{k,j}(x)
R_{j,d}
\left(
\sum_{l=-\infty}^{\infty}c_l\cdot l^d 
\right)
\beta_{\langle l \rangle}(x)
\right|
\\
&\leq&
\left\{
\left|
\frac{1}{P(x)^k}
\right|
\sum_{j=1}^k 
\sum_{d=0}^j
|Q_{k,j}(x)|
|R_{j,d}|M_d
\right\}
|\beta_{\langle 0 \rangle}(x)| \to 0 \quad (x \to 0+0)
\end{eqnarray*}
because the last $\{\cdots\}$ is a rational function 
when $x$ is close to $0$. 
Also this computation shows the validity  of the first equality. 
\hspace{\fill}$\square$

\begin{rem}{\rm \quad 
For the equation (II), 
the smoothness and flatness of 
$\beta\log \beta_{\langle 0 \rangle}$ 
for a solution $\beta$ to (I) 
follow from more or less the same 
arguments, because 
$\log \beta_{\langle 0 \rangle} 
= (R + i\theta_0) \int_{x_0}^x\frac{1}{P(u)}du$.   
}\end{rem}

\section{\large Unified proof for Case (2) and Case (3)}
The proof of Theorem \ref{Main} given in this section 
relies on the theory of center manifolds and 
the idea of graph transformation.  
For this theory, refer to a nice book by Shub 
\cite{Sh}, 
in particular, Appendix III to Chapter 5. 
First we review the center manifold theorem in 
a form which is suitable in and focused to our context. 

\begin{thm}\label{center}{\rm(Theorem III. 2, 
\cite{Sh}, 
modified) \quad 
Let $T : E \to E$ be a continuous linear endomorphism 
on a Banach space $E$ with a $T$-invariant decomposition 
$E=E_1\oplus E_2$ into closed subspaces.   
For the restrictions $T_i : E_i \to E_i$ ($i=1,\, 2$) 
of $T$ we assume that 
$T_1$ is an isomorphism and 
there exist positive constants $0<\mu^* < \lambda^*$ 
satisfying the following conditions. 
\begin{eqnarray*}
\Vert T_1({\bf v})\Vert > \mu^*\Vert {\bf v}\Vert
\quad
\mathrm{for\,\,all}\,\, {\bf v}\ne 0 \in E_1 \, , 
\\
\Vert T_2({\bf v})\Vert < \lambda^*\Vert {\bf v}\Vert
\quad
\mathrm{for\,\,all}\,\, {\bf v}\ne 0 \in E_2 \, . 
\end{eqnarray*}
Then there exists a real number $\varepsilon^*>0$ 
such that if a $C^r$-map $\Phi : E \to E$ 
($r\geq 1$) satisfies $\Phi(0)=0$ and 
$\mathit{Lip}(\Phi - T)<\varepsilon^*$, 
then we have the followings. 
\vspace{5pt}
\\
1)\quad The set 
$\,\,\displaystyle W_1= \cap_{n\geq 0} \Phi^n(S_1)\,\,$ 
where 
$\,S_1=\{({\bf v}_1, {\bf v}_2)\in E_1\times E_2 \, ; \,
\Vert {\bf v}_1\Vert \geq \Vert {\bf v}_2\Vert \}\,$  
is the graph of a $C^1$-map $\,g:E_1 \to E_2\,$ 
with 
$\mathit{Lip(g)}\leq 1$ and is invariant by $\Phi$, 
namely, 
$\Phi(W_1)=W_1\,$. 
\vspace{3pt}
\\
2)\quad If $\,\lambda^*<(\mu^*)^r\,$ holds, then the map $\,g\,$ 
is of $\,C^r$.    
}
\end{thm}
In this theorem $\mathit{Lip}$ denotes the Lipschitz constant 
of a Lipschitz map, \IE 
$\displaystyle \mathit{Lip}(f) 
= \sup\{\Vert f({\bf v_2})-f({\bf v_1})\Vert / 
\Vert{\bf v_2}-{\bf v_1} \Vert 
\, ; \, {\bf v_2}\ne {\bf v_1}, \,\,
{\bf v_1}, {\bf v_2} \in E
\}\,
$. 

In order prove Theorem \ref{Main} concerning Equation (I), 
we take $E_1=\R$, $E_2=\C$, $T_1=\mathit{id}_\R$, and 
$T_2$ is a scalar multiplication by $\lambda^{-1}\,$.  
Thereofre the real number $\varepsilon^*$ 
in the theorem is determined by $\lambda$. 

The very virtue of this theorem is that 
higher order regularities are assured 
only by estimates 
on 1-jets.  
\vspace{3pt}

Let us explain a rough idea before getting into 
the details. 
As $\Phi$, the map 
$(x, z)\mapsto (\tilde\varphi^{-1}(x), \lambda^{-1}z)$ 
or its modification  will be taken. Here 
$\tilde\varphi(x)=\varphi(x)$ for $x\geq 0$.  
If we apply this theorem by taking 
$\Phi(x, z)=(\tilde\varphi^{-1}(x), \lambda^{-1}z)$, 
we just obtain $g(x)\equiv 0$ and nothing more.  

The basic strategy is, not exactly but roughly; 
for any $\beta \in \mathcal Z$ 
and for any $\varepsilon>0\,$, 
we look for an appropriate $\Phi$ with 
$\mathit{Lip}(\Phi - T)<\varepsilon$, 
so that the resultant $g$ coincides with $\tilde\beta$ 
for $x<\delta$ for some $\delta>0$. 
Here $\tilde\beta$ is an extension of $\beta$ 
to $\R$ by taking $\tilde\beta|_{(-\infty, 0]}\equiv 0$.   
Before these arguments, 
we need to take appropriate modifications of $\varphi$, 
and for $\beta$ we choose $\Phi$ in a suitable way.  
\vspace{5pt}

Let us start the proof of the theorem for Equation (I).  
We fix a smooth function $h\in C^\infty(\R;[0,1])$ 
satisfying 
$$
h(x)\equiv 0\,\,\, \mathrm{on} \,\,(-\infty, 
1/3]
\quad 
\mathrm{and} 
\quad
h(x)\equiv 1\,\,\, \mathrm{on} \,\,[
2/3,\infty)\, .
$$

Now take $\varphi\in\Diff^\infty([0,\infty))$ 
in either of Case (2) or (3) and  
$\lambda\in \C\,$ as well.  
Take any extension of $\varphi$ in $\Diff^\infty(\R)$.  
For abuse of notation, 
it is denoted by $\varphi$ again. 
In Case (3), of course we can take $\varphi$ 
so as to be the identity on the negative side 
$(-\infty, 0]$.

First, by the following lemma, 
we 
modify it 
away from the origin so as to be 
suitable for the center manifold theory 
while its germ is not changed.

\begin{lem}\label{modification}{\rm \quad 
For $\delta>0$ define $\tilde\varphi_\delta$ as follows. 
\begin{eqnarray*} 
\tilde\varphi_\delta(x)&=&h\left(\frac{-x}{\delta}\right)\cdot x 
+ \left(1-h\left(\frac{-x}{\delta}\right)\right) \cdot \varphi(x) 
\quad \mathrm{for}\,\, x\leq 0 \, ,
\\
\tilde\varphi_\delta(x)&=&h\left(\frac{x}{\delta}\right)\cdot 
(x + \delta^2) 
+ \left(1-h\left(\frac{x}{\delta}\right)\right) \cdot \varphi(x) 
\quad \mathrm{for}\,\, x\geq 0 \, .
\end{eqnarray*} 
Then  we have
$$
\lim_{\delta\to 0+0} 
\mathit{Lip}(\tilde\varphi_\delta - \mathit{id}_\R) = 0\, .
$$
In particular, we have the followings. 
\vspace{5pt}
\\
1)\quad For small enough $\,\delta$, 
$\tilde\varphi_\delta\,$ is in $\,\Diff^\infty(\R)\,$  
and expanding on $[0,\infty)$. 
\vspace{3pt}
\\
2)\quad The germ of  $\tilde\varphi_\delta$ around $x=0$ 
is the same as that of $\varphi\,$.  
\vspace{3pt}
\\
3) \quad
$\tilde\varphi_\delta|_{(-\infty, -\delta]}
=\mathit{id}_{(-\infty, -\delta]}\,\,$  
and 
$\,\,\tilde\varphi_\delta|_{[\delta, \infty)}
=\mathit{id}_{[\delta, \infty)} + \delta^2\,$.   
\vspace{3pt}
\\
4) \quad $\ds
\lim_{\delta\to 0+0} 
\mathit{Lip}(\tilde\varphi_\delta^{-1} - \mathit{id}_\R) = 0\, .
$
}
\end{lem}

\noindent 
{\it Proof. }\quad From a direct computation using 
$\,\varphi'(0)=1$, 
the uniform convergence 
$\,\tilde\varphi_\delta'(x) \to 0\,$ 
when $\,\delta\to 0+0\,$ is 
easily obtained. 
Then, because the support of $\,\,\tilde\varphi_\delta'- 1\,\,$ 
is contained in $\,[-\delta, \delta]\,$ 
and $\,\tilde\varphi_\delta\,$ converges to 
$\,\mathit{id}_\R\,$  
uniformly, we obtain the above estimate for the Lipschitz 
constant. The statements 1) - 4) follow naturally. 
\qed

For each 
$\,\tilde\varphi_\delta|_{[0, \infty)}
\in \Diff^\infty([0, \infty))\,$,  
take any solution 
$\beta \in {\mathcal Z}_{\lambda, {\tilde\varphi_\delta}}\,$  
and the extension $\,\tilde\beta\,$ to $\,\R\,$ 
as explained above. 
Our objective is to prove that 
$\tilde\beta$ is smooth on $\,\R\,$. 
At least for $\tilde\beta'(0)$, 
not only we see it eaily 
but we need it for our proof. 
This fact is true even for the case 
$\varphi'(0)> 1$ as far as 
$\varphi'(0)<|\lambda|\,$.  
\begin{prop}\label{1jet}{\rm \quad 
$\tilde\beta$ is of $C^1\!$, namely, 
$\ds \lim_{x\to 0+0}\beta(x)=\lim_{x\to 0+0}\beta'(x)=0\,$ 
holds.  
}
\end{prop}
\noindent{Proof.}\quad 
Only $\ds \lim_{x\to 0+0}\beta'(x)=0\,$ 
is verified.   
From Equation (I), we have
$$
\beta'(\varphi(x))=(\varphi'(x))^{-1}\lambda\beta'(x)\, .
$$
From the condition there exist 
$\,x_1\!>\!0\,$ and $\,\nu\!>\!1\,$ 
such that 
$\,\varphi^{-1}|\lambda|\!<\!\nu\,$ holds 
on $\,[0,x_1]\,$.  
Take $M=\max|\beta'(x)|$ on the fundamental domain 
$[\varphi^{-1}(x_1), x_1]$ 
of the action of $\varphi$ on $(0, \infty)\,$,  
we see that when $\,x\,$ approaches to $0$ in $(0, x_1)$, 
each time it passes through a smaller fundamental domain, 
$|\beta'(x)|$ shrinks by $\nu^{-1}\,$. \qed

Next step is to look for a suitable 
$\Phi : \R\times \C \to \R\times \C \,$.  
First put 
$\Phi_0(x,z)=({\tilde\varphi_\delta}^{-1}(x), \lambda^{-1} z)$.  
The graph of $\,\tilde\beta\,$ 
is invariant under $\,\Phi_0\,$,  
while the only invariant one contained in $W_1$ 
in the center manifold theorem is the real axis 
$\R \times \{0\}\, $, because any non-trivial 
solution $\beta$ grows exponentially.  
To avoid this inconvenience, 
consider a diffeomorphism 
$$H_c(x,\, z) = \left(x, \,
z + h\left(\frac{x}{c}\right)\tilde\beta(x)\right)$$ 
of $\,\R\times \C\,$ depending on 
the parameter $\,c>0\,$.   
$H_c$ is the identity on $\{x\leq c/3\}\,$.  
Then by $H_c$ we take the conjugate 
$$\Phi_c = H_c^{-1}\circ \Phi_0 \circ H_c\, .$$

\begin{lem}\label{conjugation}
{\rm \quad For any $\,\varphi$, 
small enough $\,\delta$, 
and $\beta$, the followings hold.  
\vspace{5pt}
\\
1)\quad $\ds \lim_{c\to 0+0}
\mathit{Lip}(\Phi_c - \Phi_0)=0\,$. 
\vspace{3pt}
\\
2)\quad The graph of 
$\,\,(1-h(\frac{x}{c}))\tilde\beta(x)\,\,$ is 
invariant under $\,\Phi_c\,$.  
}
\end{lem}

\noindent
{\it Proof.}\quad 
Even though  $\,H\,$ is exponentially 
away from the identity according to $\beta$ growing 
when $\,x\to \infty\,$,   
thanks to the fact that 
$\,H\,$ and $\,\Phi_0\,$ commute 
to each other on $\,\{x\geq \tilde\varphi_\delta(c)\}\,$, 
the result of conjugation does not go away.  
Precisely the lemma is proved by direct computations as follows. 

From the definition and Equation (I) we have 
\begin{eqnarray*}
H_c^{-1}\circ\,\Phi_c\!\!\!\!&\circ&\!\!\!\! H_c(x, z) 
\,- \, \Phi_0(x, z) 
\\
&=&
\left(0\,,\,\,
 \lambda^{-1}h\left(\frac{x}{c} \right)\tilde\beta(x)
- h\left(\frac{\tilde\varphi_\delta^{-1}(x)}{c} \right)
\tilde\beta(\tilde\varphi_\delta^{-1}(x))
\right)
\\
&=& 
\left(0\,,\,\,
\lambda^{-1}\tilde\beta(x)
\left(
h\left(\frac{x}{c} \right) 
- h\left(\frac{\tilde\varphi_\delta^{-1}(x)}{c}
   \right)
\right)\right)\, .
\end{eqnarray*}
Therefore in order to conclude 
$\ds \lim_{c\to 0+0}
\mathit{Lip}(\Phi_c - \Phi_0)=0$, 
it is enough to show the uniform convergence 
of the derivative of the second component with respect to $x$ 
to $0$ when c tends to $0$. 

The estimates concerning $c\to 0$ 
which appear below are uniform in $x$.  
Let us make this point clearer. 
The second component of the above has  the support 
contained in $[0, \tilde\varphi_\delta(c) ]$ 
as a function on $x$.  
As we assumed that 
$\ds \lim_{x\to 0+0}\tilde\varphi_\delta(x) \to 1$, 
taking $c>0$ small enough, 
we can also assume that 
$\tilde\varphi_\delta(c) \leq 2c\,$ and 
it is enough to verify the estimates on $[0, 2c]\,$.
Now for example, 
as we remarked in the above proposition we know  
$|\tilde\beta(x)|=o(|x|)$ and hence 
we have 
$\max\{|\beta(x)|\, ; \, 0\leq x \leq 2c 
\}
= o(c)\,$.  

Now we show the derivative of the second component 
with respect to $x$ uniformly converges to 0 
(namely $o(1)$) when $c\to 0$.  
We can forget about $\lambda^{-1}$ 
because it is just a constant. 
\begin{eqnarray*}
\frac{d}{dx}
\left\{ 
\tilde\beta(x) 
\left(
h\left(\frac{x}{c} \right) 
- h\left(\frac{\tilde\varphi_\delta^{-1}(x)}{c}\right)
\right)
\right\}
\qquad
\qquad
\qquad
\qquad
\qquad
\quad
{}
\\
=\,
\tilde\beta'(x)
\left(
h\left(\frac{x}{c} \right) 
- h\left(\frac{\tilde\varphi_\delta^{-1}(x)}{c}\right)
\right)
\qquad
\qquad
\qquad
{}
\\ 
\qquad
\qquad
+
\,\,
\frac{\tilde\beta(x)}{c}
\left(
h'\left(\frac{x}{c} \right) 
- 
(\tilde\varphi_\delta^{-1})'(x)
h'\left(\frac{\tilde\varphi_\delta^{-1}(x)}{c}\right)
\right)
\, .
\end{eqnarray*}
The first term is of $o(1)$  because $\tilde\beta'(x)=o(1)$ 
and 
$|h(\frac{x}{c})
- h(\frac{\tilde\varphi_\delta^{-1}(x)}{c})|\leq 1\,$.  
As to the second term, 
$h'=O(1)$,  $\tilde\beta'=O(1)$, 
and $\tilde\beta(x)=o(c)$ as remarked above. 
Therefore 1) is proved. 
\vspace{5pt}

A direct computation verifies 2).  
It is also understood from the arrangement 
of $\Phi_c$, 
because $H_c$ sends the graph of  
$\,\,(1-h(\frac{x}{c}))\tilde\beta(x)\,\,$ 
to that of $\tilde\beta(x)$ 
which is invariant under $\Phi_0$. 
\qed

Now we are ready to apply the center manifold theorem 
to proof of Theorem \ref{Main} 1). 

Recall that our 
$\,T\,$ was fixed as $\,T(x, z)=(x, \lambda z)\,$.  
From out settings, 
we can take such $\mu^*$ and $\lambda^*$ that 
\vspace{3pt}
\\
\quad \bul 
$\,\,\mu^*$ is as close to 1 as we want 
as far as $\mu^*<1$ is satisfied  
\\
and
\\
\quad \bul 
$\,\,\lambda^*$ is as close to $\lambda^{-1}$ as we want 
as far as $\mu^* > \lambda^{-1}$ is satisfied.  
\vspace{3pt}
\\
Therefore for an arbitrary fixed $\,r\in \N\,$ 
we take 
$\mu^*$ and $\lambda^*$ 
so that 
$\lambda^*<(\mu^*)^r$ 
is also satisfied.  
From these choices 
$\,\varepsilon^*\,$ 
in the center manifold theorem is also fixed.  
Then by Lemma \ref{modification} 
we can find $\delta>0$ so that 
$\ds\,\mathit{Lip}(\tilde\varphi_\delta - \mathit{id}_\R) 
< \varepsilon^*/2 \,$.  
This means $\,\mathit{Lip}(\Phi_0 - T) 
< \varepsilon^*/2\,$ because 
$\ds\,\mathit{Lip}(\Phi_0 - T) = 
\mathit{Lip}(\tilde\varphi_\delta - \mathit{id}_\R)\, 
$.  
For this $\,\delta$, by Lemma \ref{conjugation} 
we can find $\,c>0\,$ so that 
$\ds\,\mathit{Lip}(\Phi_c - \Phi_0) < \varepsilon^*/2\,$.  
Therefore the center manifold theorem is applicable 
to $T$ and to our $\Phi_c$ as $\Phi$ in the theorem.  

The second statement of Lemma \ref{conjugation} tells 
that the graph of 
$\,\,(1-h(\frac{x}{c}))\tilde\beta(x)\,\,$ 
is invariant by $\,\Phi_c\,$.  
From the arrangement it is also clear that 
the graph is contained in the sector $S_1$ 
in the center manifold theorem.  
Therefore the graph is nothing but $W_1$ 
in the center manifold theorem and $g$ 
in the theorem turns out to be 
$\,\,(1-h(\frac{x}{c}))\tilde\beta(x)\,\,$ 
in our case.  
Therefore it is conluded that this function 
is of $C^r\,$ and so is $\tilde\beta$. 

We are free to improve the choice of 
$\mu^*$ and $\lambda^*$ to obtain 
another arbitrary $r\in \N\,$.  
As a conclusion, $\tilde\beta$ 
is of $C^\infty$.  
This implies nothing but 
the fact that 
$\beta$ is smooth and flat at $x=0\,$   
and completes the proof. 
\vspace{5pt}

It is easy to arrange the proof 
for Equation (II').  
The space $\,E_2\,$ is now taken to be $\,\C^M\,$ 
and the operator 
$\,T\,$ to be $\,A^{-1}\,$. 
We should remark here that 
by change of basis, $\,A\,$ can be conjugate to one 
which is arbitrarily close to 
$\,\lambda\cdot E\,$ where $E$ denotes 
the identity matrix. 
This enables us to choose $\mu^*$ and $\lambda^*$ 
in the same say as in the above proof.

\begin{rem}{\rm
\quad Instead of using the center manifold theorem, 
we can also arrange the proof so as to rely on the 
$C^r$ section theorem due to Hirsch-Pugh-Shub 
(cf. \cite{Sh}), 
which even proves the center manifold theorem.   
We also need the conjugation by $H_c$ to obtain $\Phi_c$ 
and then we may look at the space of functions 
which are supported on a large enough interval $[-R, R]$. 
}
\end{rem}

\begin{minipage}{70mm}
{Tomohiro HORIUCHI
\\
{\small Department of Mathematics, 
\\
Chuo University
\\
1-13-27 Kasuga Bunkyo-ku, Tokyo, 
\\ 112-8551, Japan
\\
horiuchi@gug.math.chuo-u.ac.jp
}}
\end{minipage}
\begin{minipage}{65mm}
{Yoshihiko MITSUMATSU 
\\
{\small Department of Mathematics, 
\\
Chuo University
\\
1-13-27 Kasuga Bunkyo-ku, Tokyo, 
\\
112-8551, Japan
\\
yoshi@math.chuo-u.ac.jp
}}
\end{minipage}

\end{document}